\documentclass[proc]{edpsmath}

\usepackage{format-rapport}
\usepackage[disable]{todonotes} % use the option [disable] for getting rid of the notes
\usepackage[ruled,linesnumbered]{algorithm2e}
\usepackage{ulem}
\usepackage{url}
\usepackage{cleveref}
\usepackage{mathabx}

\graphicspath{{figures/}}

\newcommand{\pf}{\texttt{p4est}}

\newcommand{\tpd}[2]{\partial_{#2} #1}

\newcommand{\vect}[1]{\mathbf{#1}}
\newcommand{\grad}[1]{\nabla p}

\newcommand{\norm}[1]{\Vert #1 \Vert}

\newcommand{\baseDeux}[1]{\overline{#1}^2}

\usetikzlibrary{patterns}

\begin{document}

\title{Experimenting with the p4est library for AMR simulations of two-phase flows.}
\author{Florence~Drui}\address{Laboratoire EM2C, CNRS UPR 288 et École Centrale Paris, Grande Voie des Vignes, 92295 Châtenay-Malabry Cédex}\secondaddress{Maison de la Simulation, CNRS USR 1441, Bat. 565 - Digitéo, CEA Saclay, 91191 Gif-sur-Yvette}
\author{Alexandru~Fikl}\sameaddress{2}
\author{Pierre~Kestener}\sameaddress{2}
\author{Samuel~Kokh}\sameaddress{2}\secondaddress{CEA/DEN/DANS/DM2S/STMF, CEA Saclay, 91191 Gif-sur-Yvette}
\author{Adam~Larat}\sameaddress{1}\secondaddress{Fédération de Mathématiques de l'École Centrale Paris, CNRS FR 3487, Grande Voie des Vignes, 92295 Châtenay-Malabry Cédex}
\author{Vincent~Le~Chenadec}\sameaddress{1}
\author{Marc~Massot}\sameaddress{1,4}

\begin{abstract}
Many physical problems involve spatial and temporal inhomogeneities that require a very fine discretization
in order to be accurately simulated. Using an adaptive mesh, a high level of resolution  is 
used in the appropriate areas while keeping a coarse mesh elsewhere. This idea allows 
to save time and computations, but
represents a challenge for distributed-memory environments. The MARS project (for Multiphase Adaptative 
Refinement Solver) intends to assess the parallel library \pf~for adaptive mesh, in a case of a finite volume 
scheme applied to two-phase flows. Besides testing the library's performances, particularly for load 
balancing, its user-friendliness in use and implementation are also exhibited here. 
First promising 3D simulations are even presented.
\end{abstract}

\begin{resume}
De nombreux problèmes physiques mettent en jeu des inhomogénéités spatiales et 
temporelles, qui pour être simulées correctement 
requi\`erent une discrétisation tr\`es fine.
Utiliser un maillage adaptatif pour obtenir ce niveau de  résolution  dans les zones où elle est 
requise, et garder un maillage grossier en-dehors, permet d'intéressantes économies en temps 
et ressources de calcul, mais repr\'esente un d\'efi pour le calcul distribu\'e. Le projet MARS 
(\textit{Multiphase Adaptative Refinement Solver}) a pour objectif d'\'evaluer la librairie parall\`ele 
de maillage adaptatif  \pf , appliquée à un schéma de volumes finis pour un mod\`ele 
bifluide d'\'ecoulement diphasique. Outre les performances de cette librairie,
en particulier en terme d'\'equilibrage de charge, 
sa facilité d'utilisation et d'implémentation sont mises en avant. Des premières simulations 3D 
prometteuses sont présentées.
\end{resume}

\maketitle
\setcounter{tocdepth}{1}

%/\/\/\/\/\/\/\/\/\/\/\/\/\/\/\/\/\/\/\/\/\/\/\/\/\/\/\/\/\/\/\/\/\/\/\/\/\/\/\/\/\/\/\/\/\/\/\/\/\/\/\/\/\/\/\/\/\/\/\/\/\/\/\/\/\/\/\/\/\/\/\/\/\/\/\/\/\/\
%===================================================================
%/\/\/\/\/\/\/\/\/\/\/\/\/\/\/\/\/\/\/\/\/\/\/\/\/\/\/\/\/\/\/\/\/\/\/\/\/\/\/\/\/\/\/\/\/\/\/\/\/\/\/\/\/\/\/\/\/\/\/\/\/\/\/\/\/\/\/\/\/\/\/\/\/\/\/\/\/\/\
\section{Introduction}

Adaptive Mesh Refinement (AMR) methods have been developed to  solve
problems dealing with phenomena appearing at multiple and very different spatial  and temporal scales. 
It is especially useful in the resolution of the dynamics of localized fronts or interfaces in plasma 
physics, reactive and complex flows \cite{duarte:tel-00667857}.
Combustion problems usually involve a very thin and localized flame front coupled to the 
hydrodynamics of the flow \cite{duarte:hal-00727442}. 
In this project, we are more specifically looking at problems 
of diffuse interface modeling for two-phase flows, where a good precision is
needed for describing the dynamics of the interface between the two phases.

One of the first comprehensive descriptions of AMR was given in \cite{berger84},
with an application to hyperbolic partial differential equations. This paper 
was followed by an extension of the method which accounts for the
presence of shocks and greatly simplified the AMR-related algorithms\cite{berger89}. 
Since then, AMR concepts have been implemented 
in dedicated application codes such as \texttt{RAMSES}~\cite{teyssier02} in astrophysics,
or \texttt{Gerris}~\cite{popinet03} for fluid and two-phase flow studies, among others.
They mostly follow the ideas of the Fully Threaded Tree of Khokhlov~\cite{khokhlov}. 
These tree-based AMR techniques are close, in terms of implementation to 
the fully adaptive multiresolution scheme (MR) introduced in \cite{Cohen03,muller} from the ideas 
of Harten \cite{harten95} and used 
in the \texttt{MBARETE, z-code} codes \cite{duarte:hal-00727442,sechelles2015} for various applications.

Unfortunately, these dedicated AMR or MR codes often lead to complex software designs 
due to the methods employed and 
suffer heavy difficulties in the development of new applications since the numerical 
method and the AMR technique are closely entangled. In particular, the problem of 
domain decomposition and load balancing for parallel computing 
in both  shared and distributed memory architectures
 is very delicate and 
necessitates costly complex methods coming from graphs partitioning theory 
\cite{brix2009parallelisation,sechelles2015}. 
Recently, several exceptions offering generic "\textit{hands-off}" AMR frameworks have appeared 
such as \texttt{CHOMBO}~\cite{chombodesign}
or \pf~\cite{burstedde11}. These can be used with any solver or numerical scheme and 
do not necessarily rely on the Fully Threaded Tree ideas. 

Many other AMR codes exist and we suggest the interested reader to look at 
the survey article \cite{Dubey20143217} for a large review of several such frameworks.

In the context of the MARS project, we have created an interface between a 
two-phase flow finite volume solver and a dedicated cell-based AMR library: \pf~\cite{burstedde11}.
The objectives are to test the user-friendliness and the performances of the library
in the context of HPC and to reveal its advantages and disadvantages compared with simple,
non-adaptive algorithms.

The first part of this proceeding describes the principles of AMR methods
and, more specifically, the way \pf~works. In the second part, we present the
two-fluid homogeneous equilibrium model, derived from \cite{ONERA-chante}, that
involves a system of three conservative equations (in 1D) for an isothermal
system of two fluids. The discretization of these PDEs is performed through
finite volume techniques and a Godunov-like scheme: the Riemann problems at
cells interfaces are solved approximately thanks to Suliciu's relaxation method
(see \cite{bouchut, suliciu}). Second order approximations using a MUSCL-Hancock scheme
(see \cite{leer84}) were also tested.

Several test cases are performed in order to evaluate different aspects of
\pf~library and of the numerical methods. They are presented in
the last part of this paper, with dedicated implementation details. 
The two-phase model is finally tested in some realistic
2D and 3D configurations.

%/\/\/\/\/\/\/\/\/\/\/\/\/\/\/\/\/\/\/\/\/\/\/\/\/\/\/\/\/\/\/\/\/\/\/\/\/\/\/\/\/\/\/\/\/\/\/\/\/\/\/\/\/\/\/\/\/\/\/\/\/\/\/\/\/\/\/\/\/\/\/\/\/\/\/\/\/\/\
%===================================================================
%/\/\/\/\/\/\/\/\/\/\/\/\/\/\/\/\/\/\/\/\/\/\/\/\/\/\/\/\/\/\/\/\/\/\/\/\/\/\/\/\/\/\/\/\/\/\/\/\/\/\/\/\/\/\/\/\/\/\/\/\/\/\/\/\/\/\/\/\/\/\/\/\/\/\/\/\/\/\

\section{\pf~Library: Description and Specificities}

%===================================================================

\subsection{Presentation of AMR Techniques}

\begin{figure}%[!ht]
  \begin{center}
    \null \hfill
      \subfloat[Example of overlapped grids illustrating block-based AMR.\label{fig-amrblock}]{\begin{tikzpicture}[scale=0.5]

\draw[white, fill=white] (-9,0.5) rectangle (9,0);

\begin{scope}[yslant=0.5, xslant=-1]
\draw [black, fill=blue!20, thick] (0, 0) grid (5, 5) rectangle (0, 0);
\end{scope}

\begin{scope}[yshift=2.2cm, yslant=0.5, xslant=-1]
\draw [black, thick, fill=red!20, step=0.5] (0, 1.99) grid (5, 5) rectangle (0, 2);
\end{scope}

\begin{scope}[yshift=4.2cm, yslant=0.5, xslant=-1]
\draw [black, thick, fill=green!20, step=0.25] (2.99, 2.99) grid (5, 5) rectangle (3, 3);
\end{scope}

\draw [thick, red] (-5, 4.7) -- (-5, 2.5);
\draw [thick, red] (-2, 3.2) -- (-2, 1);
\draw [thick, red] (3, 5.7) -- (3, 3.5);

\draw [thick, red] (-2, 8.2) -- (-2, 6.2);
\draw [thick, red] (0, 7.2) -- (0, 5.2);
\draw [thick, red] (2, 8.2) -- (2, 6.2);
\end{tikzpicture}}
    \hfill
      \subfloat[Example of cell refinement process illustrating cell-based AMR.\label{fig-amrcell}]{\begin{tikzpicture}[scale=0.5]

\draw[white, fill=white] (-7,0.5) rectangle (7,0);

\begin{scope}[yshift=-6.2cm, yslant=0.5, xslant=-1]
\draw [black, thick, fill=orange!20, step=0.75] (0, 0) grid (1.5, 1.5) rectangle (0, 0);
\draw [black, thick, fill=orange!20, step=0.375] (0, 0) grid (1.5, 0.75);
\draw [black, thick, fill=orange!20, step=0.75] (1.5, 0) grid (3, 1.5) rectangle (1.5, 0);
\draw [black, thick, fill=orange!20, step=1.5] (0, 1.5) grid (1.5, 3) rectangle (0, 1.5);
\draw [black, thick, fill=orange!20, step=0.75] (1.5, 1.5) grid (3, 3) rectangle (1.5, 1.5);
\draw [black, thick, fill=orange!20, step=0.375] (1.5+0.75, 1.5) grid (3, 3);
\end{scope}

\begin{scope}[yshift=-4.2cm, yslant=0.5, xslant=-1]
\draw [black, thick, fill=green!20, step=0.75] (0, 0) grid (1.5, 1.5) rectangle (0, 0);
\draw [black, thick, fill=green!20, step=1.5] (1.5, 0) grid (3, 1.5) rectangle (1.5, 0);
\draw [black, thick, fill=green!20, step=1.5] (0, 1.5) grid (1.5, 3) rectangle (0, 1.5);
\draw [black, thick, fill=green!20, step=0.75] (1.5, 1.5) grid (3, 3) rectangle (1.5, 1.5);
\end{scope}

\begin{scope}[yshift=-2.2cm, yslant=0.5, xslant=-1]
\draw [black, thick, fill=red!20, step=1.5] (0, 0) grid (3, 3) rectangle (0, 0);
\end{scope}

\begin{scope}[yslant=0.5, xslant=-1]
\draw [black, fill=blue!20, thick, step = 3] (0, 0) grid (3, 3) rectangle (0, 0);
\draw [black,dashed] (0,0)--(-6.2,-6.2);
\draw [black,dashed] (3,0)--(-6.2+3,-6.2);
\draw [black,dashed] (0,3)--(-6.2,-6.2+3);
\end{scope}

\end{tikzpicture} }
    \hfill \null
    \caption{Illustration of two approaches to locally refined meshes.}
    %\label{fig-block-cell}
  \end{center}
\end{figure}
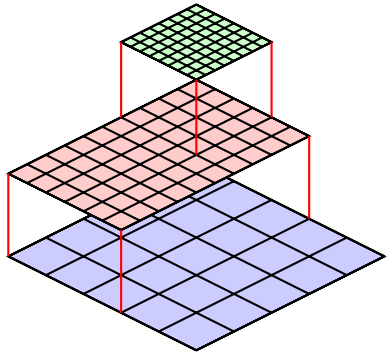
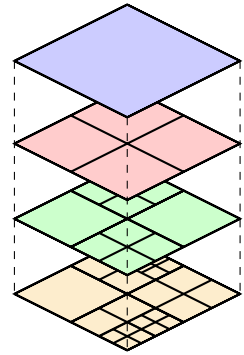

There are multiple approaches for adapting a mesh to a specific problem, among which 
block-based AMR (see Figure \ref{fig-amrblock}), cell-based AMR (Figure 
\ref{fig-amrcell}) or \textit{Wavelet-based AMR}, also called  \textit{adaptive Multi-Resolution (MR)} \cite{Cohen03,muller,koumoutsakos2011}, which can lead to error control on the solution.
Mesh-free methods, such as the \textit{Smoothed Particle Hydrodynamics} method and their multi-scale version\cite{monaghan,koumoutsakoscottet2000,koumoutsakos2005,koumoutsakos2006,koumoutsakos2007}, 
have also been
successfully employed to offer an adapted discretization. Here, we are only interested 
in mesh-based methods, and more particularly in cell-based AMR.
The cell-based method involves modifications on an initial coarse mesh (usually
a single cell representing a rectangular domain) by means of recursively dividing its
elements into multiple sub-elements with a fixed ratio. Because of the continuously
changing mesh, it is a large departure from the usual methods involving static 
discretizations. To deal with the constant modifications, cell-based AMR employs
trees to store the mesh and easily refine and coarsen specific cells. Different 
types of trees are used to store the individual refinements of each cell: 
\textbf{binary trees} for 1D domains, \textbf{quadtrees} for 2D domains, and
\textbf{octrees} for 3D (see Figure \ref{fig-p4estz}).

%===================================================================

\subsection{Challenges in cell-based AMR}

Unlike block-based AMR, where trees are only used to handle the grid hierarchy,
cell-based AMR makes heavy use of tree structures to store the mesh and modify it.
The use of this new data structure implies new difficulties in implementing
numerical methods (new integration routines, storage strategies and load
balancing techniques need to be developed). 
Indeed, during a simulation, different pieces of information about the tree structure needs to be
accessible permanently. Such a requirement raises several issues, especially
in high performance environments:
\begin{itemize}
    \item Tree storage, that can be made in linear arrays, but then raises issues of cache locality; 
    \item Tree partitioning for a better load balancing between computing processors; 
    \item Scalable algorithms; 
    \item Representation of complex geometries and not only square or rectangular domains. 
\end{itemize}
Another challenge is to include 
these tools into existing codes that need to preserve their original data structures.

Over time, various implementation choices have been made to deal with these
issues. Recently, \emph{linear tree storage}, in the form of hash-maps or linear arrays,
has been preferred to pointer-based tree representation \cite{ji,duarte:tel-00667857,sechelles2015}. These 
solutions have proved to use less memory than, for example, Fully Threaded 
Trees \cite{khokhlov}, have a better cache locality and are easier to parallelize.

\pf~is one recent example of a cell-based AMR implementation that uses linear
storage given by a space-filling curve. 
The primary usage of space filling curves in numerical simulation, 
is to provide a simple and efficient way of partitioning data for load balancing 
in distributed computing, but they can also be used for organizing data memory layout as in \pf. 
Indeed, many space filling curves have a nice property called compactness, 
which can be stated as: contiguity along the space-filing curve index 
implies, to a certain extend, contiguity in the N-dimensional space of the Cartesian mesh. 
As a consequence of the compactness property, one can expect also improved AMR 
performance due to a better cache memory usage resulting from a certain degree of
preserved locality between the computational mesh and the data memory layout.
\pf~also implements specialized refining, coarsening and iterating algorithms
for its specific choice of linear storage that have proven scalability \cite{burstedde11}.
These points are real advantages in the frame of HPC where time and space
complexity have to be handled with extreme care. 
On the top of that, the 
concept of \textbf{forest of trees}, allows to use multiple deformed but conforming
and adjacent meshes (each tree), enabling, to a certain point, to represent 
complex geometries, although this method 
does not offer the same flexibility as unstructured meshes. 

%===================================================================

\subsection{Meshing and data storage}

In \pf, the discretization of a physical space $\Omega$ is represented by multiple trees, the \emph{forest},
each tree covering a subset $\Omega_k$ of the domain, fitting its geometry. 
The trees are 
based on reference cubes $[0, 2^b]^d$, where $b$ is the maximum level of refinement and 
$d$ is the space dimension.
Each cube is sent to its corresponding subset by the one-to-one transformation
$\phi_k: [0, 2^b]^d \to \Omega_k$. 
The trees, represented by their root cell, define a \emph{macro-mesh} of the
domain, while their refined cells make up a finer \emph{micro-mesh}. 
This approach enables
the user to define complex geometries and not only square or cubic domains.
For example, let us consider a 2D annular of inner radius $R$ and outer radius $2R$.
The macro-mesh is the splitting of the annular into four four-edges fourth of an
annular, $(\Omega_k)_{k = 0, \dots , 3}$.
The corresponding one-to-one transformations are then:
$$ \varphi_k : \left\{ 
\begin{array}{ccc}
	[0, 2^b]^2 & \longrightarrow & [R, 2R] \times [0, \pi / 2] \\
	(X,Y) & \longmapsto & \left( 
			r  =  R \left( X/2^b + 1\right)  \, , \,
			\theta =  \pi/2 \left( Y/2^b + k \right)
	\right)
\end{array}
\right.
$$

The cells of the macro-mesh have to be conforming: each face (and edge in 2D) 
can only be shared by at most two trees. As each tree can have its own
spatial coordinate system, the inter-tree connectivity is static and must be
explicitly defined: this means specifying shared faces, edges and corners, 
relative orientations, etc.
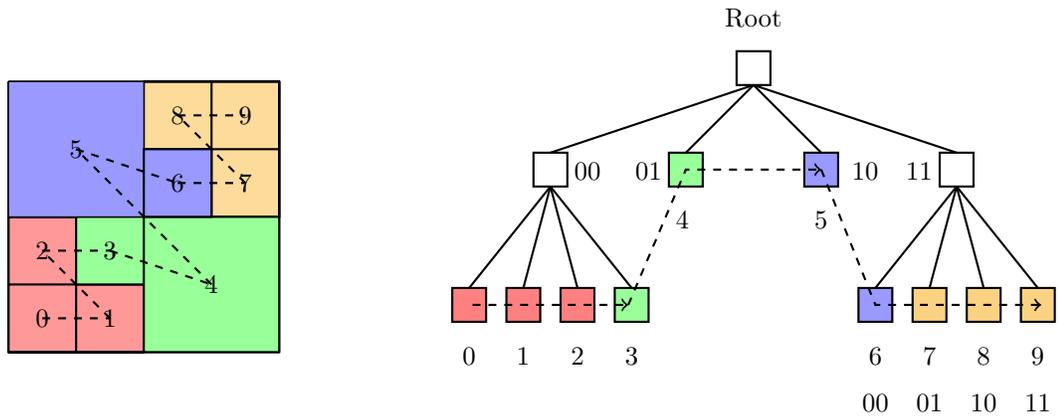
\begin{figure}%[!ht]
  \begin{center}
    \null \hfill
      \subfloat[Adaptively refined square domain. Mesh and z-order curve.]
      {\begin{tikzpicture}[scale=1.8]

\draw [white] (-1,-0.5) rectangle (3,0);

\draw [thick] (0, 0) grid (2, 2);
\draw [thick,fill=red!40,step=0.5] (0, 0) grid (1, 1) rectangle (0, 0);
\draw [thick, fill=green!40] (0.5,0.5) rectangle (1,1); 
\draw [thick,fill=orange!40,step=0.5] (1, 1) grid (2, 2) rectangle (1, 1);
\draw [thick, fill=green!40] (1, 0) rectangle (2, 1);
\draw [fill=blue!40] (0, 1) rectangle (1, 2);
\draw[thick, fill=blue!40] (1,1) rectangle (1.5,1.5);

\draw [thick, dashed] (0.25, 0.25) -- (0.75, 0.25);
\draw [thick, dashed] (0.75, 0.25) -- (0.25, 0.75);
\draw [thick, dashed] (0.25, 0.75) -- (0.75, 0.75);
\draw [thick, dashed] (0.75, 0.75) -- (1.5, 0.5);
\draw [thick, dashed] (1.5, 0.5) -- (0.5, 1.5);
\draw [thick, dashed] (0.5, 1.5) -- (1.25, 1.25);
\draw [thick, dashed] (1.25, 1.25) -- (1.75, 1.25);
\draw [thick, dashed] (1.75, 1.25) -- (1.25, 1.75);
\draw [thick, dashed] (1.25, 1.75) -- (1.75, 1.75);

\node at (0.25, 0.25) {$0$};
\node at (0.75, 0.25) {$1$};
\node at (0.25, 0.75) {$2$};
\node at (0.75, 0.75) {$3$};
\node at (1.5, 0.5) {$4$};
\node at (0.5, 1.5) {$5$};
\node at (1.25, 1.25) {$6$};
\node at (1.75, 1.25) {$7$};
\node at (1.25, 1.75) {$8$};
\node at (1.75, 1.75) {$9$};
\end{tikzpicture}
      \label{fig-zcurve}}
    \hfill
      \subfloat[The corresponding representation of the domain using a quadtree.]
      {\begin{tikzpicture}[scale=0.9]

\draw [thick] (0, 9) rectangle (0.5, 9.5);

\draw [thick] (-3, 7.5) rectangle (-2.5, 8);
\draw [thick, fill=green!40] (-1, 7.5) rectangle (-0.5, 8);
\draw [thick, fill=blue!40] (1, 7.5) rectangle (1.5, 8);
\draw [thick] (3, 7.5) rectangle (3.5, 8);

\newcommand{\xs}{-3}
\foreach \x in {\xs, \xs + 2, \xs + 4, \xs + 6} {
    \draw [thick] (\x + 0.25, 8) -- (0.25, 9);
}

\renewcommand{\xs}{-4.2}
\foreach \x in {\xs, \xs + 0.8, \xs + 1.6} {
    \draw [thick] (\x + 0.25, 6) -- (-2.75, 7.5);
    \draw [thick, fill=red!50] (\x, 5.5) rectangle (\x + 0.5, 6);
}
\draw [thick] (\xs + 2.4 + 0.25, 6) -- (-2.75, 7.5);
\draw [thick, fill=green!40] (\xs + 2.4, 5.5) rectangle (\xs + 2.9, 6);

\renewcommand{\xs}{1.8}
\foreach \x in {\xs + 0.8, \xs + 1.6, \xs + 2.4} {
    \draw [thick] (\x + 0.25, 6) -- (3.25, 7.5);
    \draw [thick, fill=orange!50] (\x, 5.5) rectangle (\x + 0.5, 6);
}
\draw [thick] (\xs + 0.25, 6) -- (3.25, 7.5);
\draw [thick, fill=blue!40] (\xs, 5.5) rectangle (\xs + 0.5, 6);

\draw[->] [thick, dashed] (-3.9, 5.75) -- (-1.6, 5.75);
\draw[->] [thick, dashed] (-1.6, 5.75) -- (-0.75, 7.75) -- (1.25, 7.75);
\draw[->] [thick, dashed] (1.25, 7.75) -- (2.05, 5.75) -- (4.5, 5.75);

\node at (0.25, 10) {Root};
\node at (-3.95, 5) {$0$};
\node at (-3.15, 5) {$1$};
\node at (-2.35, 5) {$2$};
\node at (-1.55, 5) {$3$};

\node at (-0.8, 7) {$4$};
\node at (1.25, 7) {$5$};

\node at (2.05, 5) {$6$};
\node at (2.85, 5) {$7$};
\node at (3.65, 5) {$8$};
\node at (4.45, 5) {$9$};

\node at (-2.2, 7.73) {$00$};
\node at (-1.3, 7.73) {$01$};
\node at (1.9, 7.73) {$10$};
\node at (2.7, 7.73) {$11$};
\node at (2.05, 4.3) {$00$};
\node at (2.85, 4.3) {$01$};
\node at (3.65, 4.3) {$10$};
\node at (4.45, 4.3) {$11$};
\end{tikzpicture} \label{fig-amrquadtree}}
    \hfill \null
  \end{center}
  \caption{z-order traversal of the quadrants in one tree of the forest and load partition into four processes.
      Dashed line: z-order curve. 
      Quadrant label: z-order index.
      Color: MPI processes.
  }
  \label{fig-p4estz}
\end{figure}
Each cell in the micro-mesh is then associated with 
its position in the reference cube $[0, 2^b]^d$. Therefore, each subdivision of the root node, called \textbf{octant} in 
3D (and \textbf{quadrant} in 2D) is uniquely tracked by its \emph{integer} 
spatial coordinates $(x, y, z) \in \llbracket 0, 2^b \rrbracket^3$. 
Linear storage requires a one-to-one mapping from 
the spatial coordinates $(x,y,z)$ to a linear index $m$. In \pf, this mapping is
provided by the Morton space filling curve, also called
 z-order curve. This is illustrated in Figure \ref{fig-zcurve}, where we can see
 how the z-order curve, in dashed line, covers a two-dimensional mesh.
Figure \ref{fig-amrquadtree}, illustrating the tree version of z-ordering, also shows 
an example of load balancing for four processes: each color represents a  
different process and we clearly see how the linear storage enables a simple distribution 
of the leaves of the tree.

The linear array containing all the cells is thus indexed by the Morton index $m$, constructed by
considering the binary representations of the coordinates and interwoven in the following 
way:
\beq{eq-MortonIndex}
  \baseDeux{m}_{3i + 2} = \baseDeux{z}_i,\quad \baseDeux{m}_{3i + 1} = \baseDeux{y}_i,\quad \baseDeux{m}_{3i + 0} = \baseDeux{x}_i,\quad
  \forall i \in \llbracket 0, b - 1 \rrbracket,
\eeq
where $\baseDeux{x}_i$ is the $i$-th bit of the binary representation of $x$, notation $\baseDeux{\, \cdot \,}$ 
indicating numbering in base 2. 
Using this method, the connectivity inside each tree is stored implicitly.
It enables easy finding of the direct parents, children or siblings of a given cell by
simple bit flips (details can be found in \cite{burstedde11}). However, other neighbors of
a cell from the same tree require more work to be found
because they need to be identified in the linear array where they are stored. 
This is generally achieved by iterating over the faces, edges or corners
of a given cell. 
In the case when the concerned cell stands at the boundary of the tree, neighbors in the next tree are 
found through the knowledge of its spatial coordinates and the one-to-one relation with its Morton index. However, 
one has to be careful with the possible implicit change of coordinates in the neighboring tree. 
%===================================================================

\subsection{Refining and coarsening}

\pf~creates the mesh only once, initially, and then adapt it at will
by modifying the micro-mesh. While adapting, the following
steps are usually taken:
\bitem
  \item Going through the linear array, the leaves are marked for refinement or coarsening or left unchanged, following a criterion given by the user;
  \item The refinement and coarsening is then applied on each leaf, if possible. A very important feature of the 
        adapting algorithms is that the z-order is maintained while modifying the linear array. 
        Both algorithms run in $\mathcal{O}(N_{\text{leaves}})$, but the refining
        algorithm requires additional space.
  \item \textbf{2:1 balancing} is performed: for practical reasons, the level difference between an octant 
        and each of its neighbors is at most 1 (+1 or -1), so that the neighbors of  a quadrant are at most 2 times 
        smaller or 2 times bigger; hence the 2:1 notation. Trees with such a property are also denoted "graded trees" in \cite{Cohen03,muller07,duarte:tel-00667857,duarte2012,sechelles2015}. Algorithms that perform the balancing are generally 
        among the costliest parts of AMR or Multi-resolution codes (see more details in \cite{burstedde12}). 
  \item Finally, as far as parallelization is concerned, load distribution is operated between processes by an equal division of the new array of leaves.
\eitem

%/\/\/\/\/\/\/\/\/\/\/\/\/\/\/\/\/\/\/\/\/\/\/\/\/\/\/\/\/\/\/\/\/\/\/\/\/\/\/\/\/\/\/\/\/\/\/\/\/\/\/\/\/\/\/\/\/\/\/\/\/\/\/\/\/\/\/\/\/\/\/\/\/\/\/\/\/\/\
%===================================================================
%/\/\/\/\/\/\/\/\/\/\/\/\/\/\/\/\/\/\/\/\/\/\/\/\/\/\/\/\/\/\/\/\/\/\/\/\/\/\/\/\/\/\/\/\/\/\/\/\/\/\/\/\/\/\/\/\/\/\/\/\/\/\/\/\/\/\/\/\/\/\/\/\/\/\/\/\/\/\

\section{Model and Numerical Methods}

%===================================================================

\subsection{A Two-fluid Model}\label{section- two-phase model}

We consider a two-phase flow that involves two compressible fluids $k = 1, 2$
governed by a barotropic Equation of State (EOS) $\rho_k\mapsto p_k(\rho_k)$, 
here $\rho_k$ and $p_k$ denote, respectively, the density and the pressure of
the fluid $k = 1, 2$. We make the classic assumption that $p_k^{\prime} > 0$
which enables the definition of the pure fluid sound velocities $c_k$,  $k=1,2$, by 
$c_k^2 (\rho_k)= p_k^{\prime}(\rho_k)$. The global density of the
medium is defined by:
\begin{equation}
 \rho  = \alpha\rho_1 + (1-\alpha)\rho_2,
\end{equation}
where $\alpha$ (resp. $1-\alpha$) denotes the volume fraction of fluid $1$ 
(resp. $2$). We note $Y = \rho_1 \alpha/\rho$ the mass fraction of the fluid
$1$. Following \cite{ONERA-chante,caro}, we suppose that the pressure 
$(\rho,Y)\mapsto p$ in the two-component medium is defined by imposing an
instantaneous pressure equilibrium between each fluid, namely:
\begin{equation}
 \label{eq-mechanicalEquilibrium}
 p = p_1\left(\frac{\rho Y}{\alpha}\right) = p_2\left(\frac{\rho(1-Y)}{1-\alpha}\right),
\end{equation}
for given values of $\rho$ and $Y$. The mechanical equilibrium 
relation~\eqref{eq-mechanicalEquilibrium} imposes that $\alpha$ is defined as
a function of $\rho$ and $Y$. The uniqueness of $\alpha$ verifying 
\eqref{eq-mechanicalEquilibrium} is ensured by the hypothesis
$p_k^{\prime} > 0$ and pressure is thus a function of $\rho$ and $Y$. Assumptions that ensure 
existence will be stated later, in the framework of a particular choice of equations of state, 
used in section~\ref{section- results}. 

We note $(\be_x,\be_y,\be_z)$ the canonical base of $\mathbb{R}^3$. We suppose that both fluids $k=1,2$ share the same velocity $\bu=(u_x,u_y,u_z)$, which
gives the following governing equations for the flows:
\begin{equation}\label{eq-twophasevector}
\tpd{\W}{t} + \tpd{\mathcal{F}_x(\W)}{x} +
                    \tpd{\mathcal{F}_y(\W)}{y} +
                    \tpd{\mathcal{F}_z(\W)}{z} = \S(\W),
\end{equation}
where $\W = [\rho,\rho Y,\rho u_x, \rho u_y,\rho u_z]^T$, % 
and the  fluxes $\mathcal{F}_q$ verify the rotational invariance property
$\mathcal{F}_q(\W) = R_q^{-1}\mathcal{F}_x(R_q \W)$, with 
$\mathcal{F}_x = [\rho u_x ,\rho Y u_y ,\rho u_x^2+p, \rho u_x u_y ,\rho u_x u_z ]^T$
and $R_q$ being defined as the rotation matrix:
$$
R_x = 
\left[
\begin{smallmatrix}
 1 & 0 & 0 & 0 & 0 
\\
 0 & 1 & 0 & 0 & 0 
\\
 0 & 0 & 1 & 0 & 0 
\\
 0 & 0 & 0 & 1 & 0  
\\
 0 & 0 & 0 & 0 & 1
\end{smallmatrix}
\right],
\quad
R_y = 
\left[
\begin{smallmatrix}
 1 & 0 & 0 & 0 & 0 
\\
 0 & 1 & 0 & 0 & 0 
\\
 0 & 0 & 0 & 1 & 0 
\\
 0 & 0 & 0 & 0 & 1  
\\
 0 & 0 & 1 & 0 & 0
\end{smallmatrix}
\right]
,\quad
R_z = 
\left[
\begin{smallmatrix}
 1 & 0 & 0 & 0 & 0 
\\
 0 & 1 & 0 & 0 & 0 
\\
 0 & 0 & 0 & 0 & 1 
\\
 0 & 0 & 1 & 0 & 0  
\\
 0 & 0 & 0 & 1 & 0
\end{smallmatrix}
\right].
$$
The body force term $\S$ accounts for gravity with $\S = (0,0,0,-\rho g, 0)^T$.

The system obtained by considering \eqref{eq-twophasevector} with $\S=0$ is
hyperbolic. For one-dimensional problems the resulting eigenstructure is a set
of three eigenvalues $u_x \pm c$, $u_x$ where  the sound velocity of the mixture $c^2(\rho,Y)= (\partial p /\partial \rho)_Y$ is given 
by Wood's formula~\cite{wood}:

\[
\frac{1}{(\rho c)^2}
= 
\frac{Y}{(\rho_1 c_1(\rho_1))^2} + \frac{1 - Y}{(\rho_2 c_2(\rho_2))^2}, \quad \rho_1 = \frac{\rho\,Y}{\alpha(\rho,Y)}, \quad \rho_2 = \frac{\rho(1-Y)}{(1-\alpha(\rho,Y))} 
.
\]

The fields associated with $u_x \pm c$ (resp. $u_x$) are genuinely nonlinear (resp. 
linearly degenerate). Introducing the free energy $F(\rho,Y) = \int^\rho \frac{P(r,Y)}{r^2} \text{d}r$,
the system is also equipped with a mathematical entropy inequality:
$$
\partial_t[\rho F(\rho,Y) + \rho |\bu|^2/2 + \rho g y]
+
\text{div}[(\rho F(\rho,Y) + \rho |\bu|^2/2 + P +\rho g y) \bu]
\leq 0.
$$

%===================================================================

\subsection{Finite Volume Method}

\subsubsection*{Dimensional Splitting and Finite Volume Discretization}

In order to approximate the solutions of \eqref{eq-twophasevector}, we use
a Finite Volume scheme based on a dimensional splitting strategy, namely of the
Lie splitting type. This consists, during a time step $\Delta t$, in successively 
solving one-dimensional problems, for each direction, using a discretization of 
\eqref{eq-twophasevector} through a classic 1D finite volume method. We adopt classic notations pertaining 
to unstructured meshes for describing the AMR grid: the cell $i$ is noted $K_i$ whose volume is $|K_i|$ while $|\Gamma_{ij}|$ and $\bn_{ij}$ 
are respectively the surface and the unit normal of the interface between two neighboring cells $i$ and $j$. The vector $\bn_{ij}$ 
is oriented from cell $i$ to cell $j$. We note $\mathcal{N}_q(i)$, the set of cells neighboring cell $i$ in the direction $q=x,y,z$. 
The full scheme for advancing~\eqref{eq-twophasevector} in cell $i$ from time referenced by $n$ to time $n+1$ is:
\beq{eq-update}
    \W_i^{*}  = \Box_{x}^{\Delta t} \W_i^{n}, \qquad
    \W_i^{**} = \Box_{y}^{\Delta t} \W_i^{*},  \qquad
    \W_i^{n + 1} = \Box_{z}^{\Delta t} \W_i^{**}.
\eeq
where the operator $\Box_{q}^{\Delta t}$ is defined by
\begin{equation}
\Box_{q}^{\Delta t} \W_i = \W_i 
- 
\frac{\Delta t}{|K_i|} \sum_{j \in \mathcal{N}_q(i)} |\Gamma_{ij}| (\mathbf{e}_q^T\bn_{ij}) R_q^{-1} \Phi_{ij}
,
\quad
\Phi_{ij} = \Phi(R_q\W_i,R_q\W_j),
\label{eq- numerical scheme sweep q direction}
\end{equation}
for $q=x,y,z$ and $(\W_L,\W_R)\mapsto\Phi$ being a choice of numerical flux, which has to be provided.

\subsubsection*{First order Suliciu's relaxation method}\label{section- numerical flux definition}
In order to define the numerical flux $\Phi$ we choose here the flux defined by 
the Suliciu's relaxation approach~\cite{suliciu, chalons2008}. This method belongs to the family of HLLC 
solvers~\cite{toro1994,toro} and using our notations we have
$$
\begin{aligned}
\Phi(\W_L,\W_R) = 
\frac{1}{2}
\Big[
\mathcal{F}_x(\W_L) + \mathcal{F}_x(\W_R)
&
-
\left|(u_x)_L - \frac{a}{\rho_L}\right| (\W^*_L - \W_L ) 
\\
&
- 
|u^*| ( \W^*_R - \W_L^* ) 
-
\left|(u_x)_R + \frac{a}{\rho_R}\right| (\W_R- \W_R^*)
\Big],
\end{aligned}
$$
with $u^* = \frac{(u_x)_L + (u_x)_R}{2} - \frac{1}{2a}(p_R - p_L)$, 
$1/\rho^*_L  = 1/\rho_L + \frac{u^* - (u_x)_L}{a}$, 
$1/\rho^*_R  = 1/\rho_R - \frac{u^* - (u_x)_R}{a}$, 
$Y_L^* = Y_L$,
$Y_R^* = Y_R$,
$(u_y)_L^* = (u_y)_L$,  
$(u_z)_L^* = (u_z)_L$,  
$(u_y)_R^* = (u_y)_R$,  
$(u_z)_R^* = (u_z)_R$ and $a$ is defined by
$a = \theta \max(\rho_L c_L, \rho_R c_R)$,
where $c_L$ and $c_R$ denote the sound velocity evaluated for the state $\W_L$ and $\W_R$. The parameter $\theta>1$ is a constant. This choice of $a$ complies with the subcharacteristic condition of Whitham for stability purposes (see \cite{chalons2008}).

\subsubsection*{Time step and CFL condition}

\espace

The stability of the scheme is assured at each time step by the following CFL condition:
\begin{equation}\label{eq-twophasecfl}
\Delta t \leq C \min_i \left(\frac{\Delta x_i}{\norm{\vect{u}_i} + \frac{a}{\rho_{i}}}\right),
\end{equation}
with $C\in[0,1]$. The time step is thus global over the whole mesh. Let us insist on the fact that local time stepping can be an important additional feature
in order to save computational time, which can be implemented for the resolution in time of the convective part of the system \cite{muller07,postel09}. 
Because of the framework of the original project, it is not 
considered in this contribution, but stands within our list of further improvements.

%===================================================================

\subsection{Higher-order discretization}\label{section- higher order}
We propose a simple second order extension of the Finite Volume Method presented in section~\ref{section- numerical flux definition} by 
using a classic MUSCL-Hancock strategy~\cite{leer84,berthon06} for each sweep in the direction $q=x,y,z$.

\subsubsection*{Evaluation of slopes within the cells}
We consider here only the sweep in the direction $x$: the other cases can be deduced by substituting $x$ by $y$ or $z$. Let $\Lambda$ be the change of variables $\Lambda: \W \mapsto [\rho_1 \alpha, \rho_2 (1-\alpha), u_x, u_y,u_z]=\mathbf{V}$.
For a cell $K_i$, for each $j\in\mathcal{N}_x(i)$, the set of neighbors of $K_i$ in the $x$ direction, 
we define a slope $\sigma_{ij}$ for the variations of the primitive variables in the direction $x$ 
at each interface $\Gamma_{ij}$ with $\sigma_{ij} = (\mathbf{V}_j - \mathbf{V}_i) / (\mathbf{M}_i\mathbf{M}_j \cdot \mathbf{e}_x  )
$, where 
$M_r$ is the center of the cell $K_r$ and $\mathbf{V}_r=\Lambda(\W_r)$, $r=i,j$. Then we evaluate a slope $\sigma_i$ associated with the variations along the $x$ axis 
in the vicinity of $K_i$ thanks to a simple \textit{minmod} limiting procedure that accounts for all $\sigma_{ij}$ by setting:
\[
\sigma_i 
= 
\begin{cases}
 s
 \min\{|\sigma_{ij}| , j\in\mathcal{N}_x(i) \} 
 ,\quad & \text{if all $\sigma_{ij}$ for $j\in\mathcal{N}_x(i)$ have the same sign $s=\pm 1$,}
 \\
 0
 ,\quad & \text{otherwise.}
\end{cases}
\]

\subsubsection*{Prediction step}
For a given cell $i$, the MUSCL-Hancock method involves the computation of the two
left and right predicted values  $\W_{iL}^{n+\demi}$ and $\W_{iR}^{n+\demi}$ for the conservative variables as follows:
\begin{itemize}
  \item compute $\W_{iL}^n$ and $\W_{iR}^n$ in the cell $i$ with:
      $
      \W_{iL}^n  = \Lambda^{-1} (\mathbf{V}_i^n - \frac{\sigma_i \Delta x_i}{2})
      $,
      $\W_{iR}^n  = \Lambda^{-1} (\mathbf{V}_i^n + \frac{\sigma_i \Delta x_i}{2})
      $.
  \item evaluate predicted left and right states $\W_{iL}^{n+\demi}$ and $\W_{iR}^{n+\demi}$ in the cell $i$ with :
  \begin{align}
    \W_{iL}^{n+\demi}  &= \W_{iL}^n - \demi \frac{\Delta t}{\Delta x} \left( \Flux(\W_{iR}^n) - \Flux(\W_{iL}^n) \right),& 
    \W_{iR}^{n+\demi}  &= \W_{iR}^n - \demi \frac{\Delta t}{\Delta x} \left( \Flux(\W_{iR}^n) - \Flux(\W_{iL}^n) \right), \\
  \end{align}
\end{itemize}
Let us note that the change of variables $\Lambda$ for computing $\sigma_{ij}$ ensures the positivity of $\alpha \rho_1$ and $(1-\alpha) \rho_2$.

\subsubsection*{Flux computation}
The definition of the flux $\Phi_{ij} $ at an interface $\Gamma_{ij}$ by the MUSCL-Hancock method is obtained by replacing the left and 
right values in the flux $\Phi$ by the predicted values as follows: one replaces the choice of $\Phi_{ij}$ in relation~\eqref{eq- numerical scheme sweep q direction} 
by $\Phi^\text{MH}_{ij}$ where  $\Phi^\text{MH}_{ij} = \Phi(R_x \W_{iR}^{n+\demi}, R_x \W_{jL}^{n+\demi})$ 
if $\mathbf{e}_x^T \mathbf{n}_{ij}>0$ and $\Phi_{ij} = \Phi(R_x \W_{jR}^{n+\demi}, R_x \W_{iL}^{n+\demi})$ otherwise. 

\subsubsection*{Strang splitting}
Lie splitting formulae introduce an asymptotically first order global error, whereas using Strang splitting formulae lead to a second order 
when the splitting substeps are resolved exactly 
\cite{leveque,duarte:tel-00667857}.  Second order is maintained if each substep is integrated in time with a numerical scheme at least of second order
\cite{duarte:tel-00667857}. In our case one
solves according to the $x$ then $y$ directions on a half time step, along the $z$ direction on a full time step, and again
$y$ then $x$ directions on a half time step.
The update procedure~\eqref{eq-update} is replaced by:
\beq{eq-updateStrang}
\Box_X^{\Delta t/2}
\Box_Y^{\Delta t/2}
\Box_Z^{\Delta t/2}
\Box_Z^{\Delta t/2}
\Box_Y^{\Delta t/2}
\Box_X^{\Delta t/2}
\W_i^n,
\eeq
leading to a numerical scheme that is second order in time and space.

%===================================================================

\subsection{Refinement criterion}
\label{subsec-criterion}
The definition of efficient refinement criteria is a complex 
task that depends on the physical phenomena involved in the simulation (see \textit{e.g.}\cite{duarte:tel-00667857,berger84}). We consider here only three simple heuristic criteria 
in order to test the mesh adaptation functionality of \pf{} within our finite volume framework.
The authors are aware that this part is critical in the study of an AMR technique and 
definitely  requires a deeper investigation. In the following, we briefly describe the different criteria 
we have tested so far. 

Each time the mesh-adapting algorithms are called, a given criterion $C(\W)$ is evaluated within each 
cell and compared 
to a given threshold $\xi$. If $C(\W) > \xi$, then the current cell must be refined. If all siblings of a given octant verify $C(\W) \leq \xi$, then the octant is marked for coarsening. The final configuration of the mesh is obtained by accounting for the 2:1 balance constraint. 
During coarsening, the new coarser cell contains the mean value of the to-be-removed cells. During refining,
the new cells are fed with the mean value of their parent cell, even when using second-order reconstruction.
By experiment, feeding the new cells with more accurate values has not shown substantial improvement. 
Let $b$ denote a scalar or a vector value, we note
$D(b)_i 
= 
\max 
\{ 
\frac{|b_i - b_j|}{\max(b_i,b_j)}~/~j\in \mathcal{N}_q(i), q=x,y,z 
\}
$.

In the following tabular, we have ordered the criteria by increasing sensibility. The mildest \textbf{$\alpha$-gradient}
allows to refine only at the interface. A \textbf{mixed-criterion} involving local jumps of density, velocity and pressure
with a rather high threshold $\xi$ allows to moreover captures non-linear waves and strong variations of the solutions. 
Finally, the criterion on the density only with a low threshold allows to capture all the small variations in the 
solution, even possibly the acoustic features. 

\espace

\begin{center}
\begin{tabularx}{\linewidth}{c l X }
	\hline
	\textit{Name} & \textit{{Description}} & \textit{{Use}} \\
	\hline
	{} & {} & {} \\
	\textbf{$\alpha$-gradient} & 
	$
	C(\W)_i = 
	\max \{
	  |\alpha_i - \alpha_j|~/~j\in \mathcal{N}_q(i), q=x,y,z
	  \}
	$ 
	& {\small Evolution of the interface gas/liquid.} \\
	{} & {} & {} \\
	\textbf{mixed-criterion} 
	& 
	$
	C(\W)_i =
	\max [a D(\rho)_i, b D(p)_i, c D(\vect{u})_i]
	$
	& {\small $\bigplus$ General criterion with selection of prevailing non-linear waves according to $a$, $b$, 
	$c$ weights.} \\
		\textbf{$\rho$-gradient} 
		&
		$C(\W)_i = D(\rho)_i$ & 
		{\small $\bigplus\bigplus$ Most sensible criterion. Captures all variations in the solution, even small amplitude acoustic waves. } \\
\end{tabularx}
\end{center}

%/\/\/\/\/\/\/\/\/\/\/\/\/\/\/\/\/\/\/\/\/\/\/\/\/\/\/\/\/\/\/\/\/\/\/\/\/\/\/\/\/\/\/\/\/\/\/\/\/\/\/\/\/\/\/\/\/\/\/\/\/\/\/\/\/\/\/\/\/\/\/\/\/\/\/\/\/\/\
%===================================================================
%/\/\/\/\/\/\/\/\/\/\/\/\/\/\/\/\/\/\/\/\/\/\/\/\/\/\/\/\/\/\/\/\/\/\/\/\/\/\/\/\/\/\/\/\/\/\/\/\/\/\/\/\/\/\/\/\/\/\/\/\/\/\/\/\/\/\/\/\/\/\/\/\/\/\/\/\/\/\

\section{Results}\label{section- results}
We present 2D and 3D simulations performed with the code developed during the CEMRACS 2014 
research session.
These results aim at testing several elements: the AMR functionalities of \pf{}, the computational cost
reduction thanks to the compression of the mesh and the parallel performance of \pf{}. 
We also propose simulations of gravity driven flows with the 
two-phase model of section~\ref{section- two-phase model}. 
Let us emphasize the fact that tests are early results that shall be more thoroughly 
investigated in the future.

In the sequel we shall consider that the EOS of each component $k=1,2$ is barotropic 
Stiffened Gas law of the form
\(
\rho_k \longmapsto p_k(\rho_k) = p_{k, 0} + c_k^2 (\rho_k - \rho_{k, 0}),
\)
where $p_{k, 0}, \rho_{k, 0}$ and $c_k$ are positive characteristic constants
of the fluid $k$. This choice of EOS ensures that $\alpha$ and $P$ can always be uniquely 
defined thanks to explicit formulas~\cite{ONERA-chante,caro}.

%===================================================================

\subsection{Scheme verification}\label{subsec-adv}

We consider here several tests that consist in advecting a constant velocity and constant pressure profile in a periodic 2D domain $[0,1]^2$
(all the physical dimensions will be given in SI units).

The initial condition is defined by $p(\vect{x}, 0) = 10^5$, 
$\vect{u}(\vect{x}, 0) = (1, 1, 1)^t $
 and a given initial profile of $\alpha$ defined by a function $\alpha_0(\vect{x})$. The exact solution of this problem is trivially  $\alpha(\vect{x},t)=\alpha_0(\vect{x}-t\vect{u})$ with $p$ and $\vect{u}$ kept at their initial value.

\begin{figure}
\centering
\subfloat[Illustration of $\alpha$ profile (first order in blue, MUSCL-Hancock in purple). $\Delta x_{max} = 2^{-3}$, $\Delta x_{min} = 2^{-8}$. Cut along $x=y$ line.]
{\includegraphics[width=7cm]{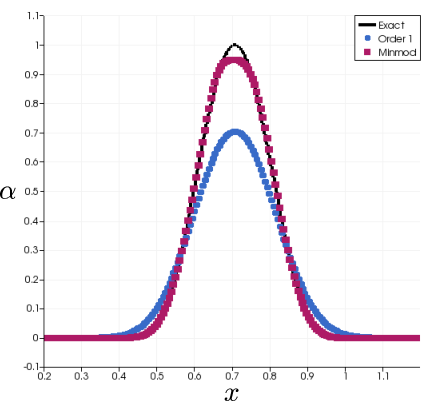}\label{gaussian-result}}
    \subfloat[Convergence rate on uniform meshes for the first-order and second-order schemes, with L$^1$ and L$^2$ norms.]{\includegraphics[width=8cm]{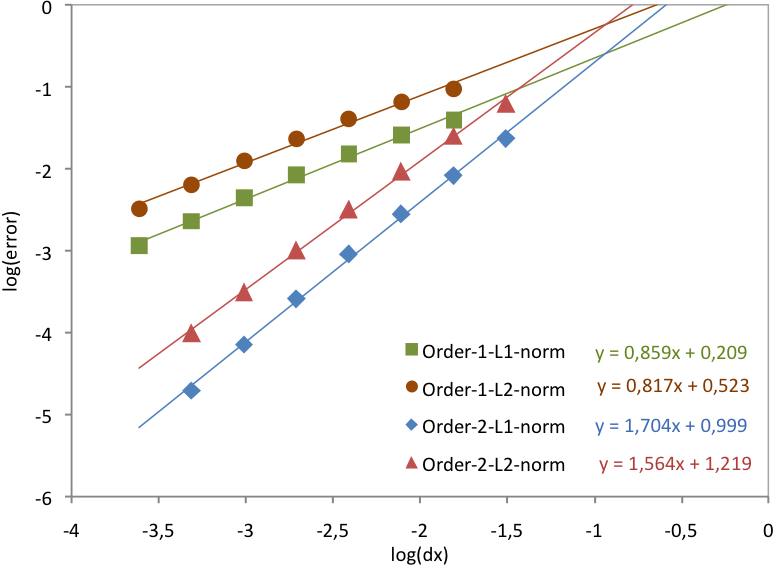}\label{convergence}}
    \caption{Advection of a smooth $\alpha$-profile, computed with first-order and second-order MUSCL-Hancock schemes after 1s of simulation.}
\end{figure}

First, we want to evaluate the behavior of the MUSCL-Hancock method with the simple slope evaluation 
described in section~\ref{section- higher order} in a AMR context. We suppose that $\alpha_0$ is given 
by a smooth profile 
\begin{equation}
\alpha_0(\vect{x}) = 
\begin{cases}
\lambda + (1-\lambda) \cdot \cos^4 \left(\displaystyle\pi \frac{|\vect{x} - \vect{x_0}|}{0.6} \right), &  \text{if } |\vect{x} - \vect{x_0}| \geq 0.3, \\
\lambda, & \text{otherwise},
\end{cases}
\label{eq- alpha profile}
\end{equation}
for $\lambda = 10^{-7}$, $\vect{x}_0=(0.5,0.5)$
and we choose to drive the AMR with the \textbf{$\rho$-gradient} criterion with a threshold value of $\xi = 5.10^{-5}$.
Figure~\ref{gaussian-result} shows the resulting $\alpha$-profile at $t=1$,
with a space step ranging from $\Delta x_{\max{}} = 2^{-3}$ to  $\Delta x_{\min{}} = 2^{-8}$. As expected, the higher-order method clearly improves the accuracy of the solution. In figure~\ref{convergence}, we verify that the
convergence rate in the $L^1$ and $L^2$ norms are compatible with the standard results
\cite{leveque}.
The proposed evaluation, involving points for coarse meshes where order reduction in the non-asymptotic regime is taking place, still gives a convergence rate of 0.8 for the first order scheme and 1.6 for our MUSCL-Hancock implementation.

%===================================================================
\subsection{Tests of parallel AMR procedure}
\begin{figure}
  \begin{center}
    \subfloat{\includegraphics[height=0.35\tw]{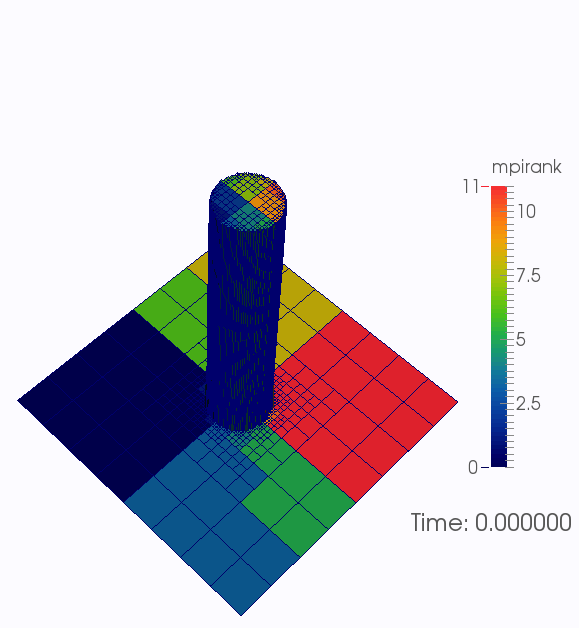}}
    \subfloat{\includegraphics[height=0.35\tw]{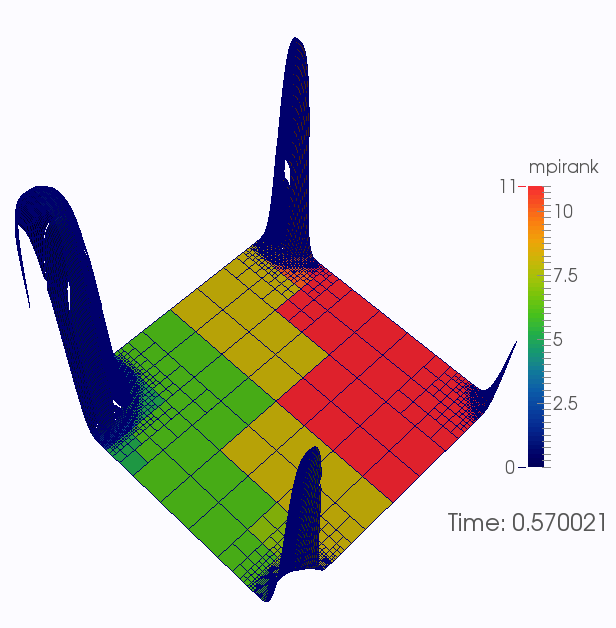}}
    \subfloat{\includegraphics[height=0.35\tw]{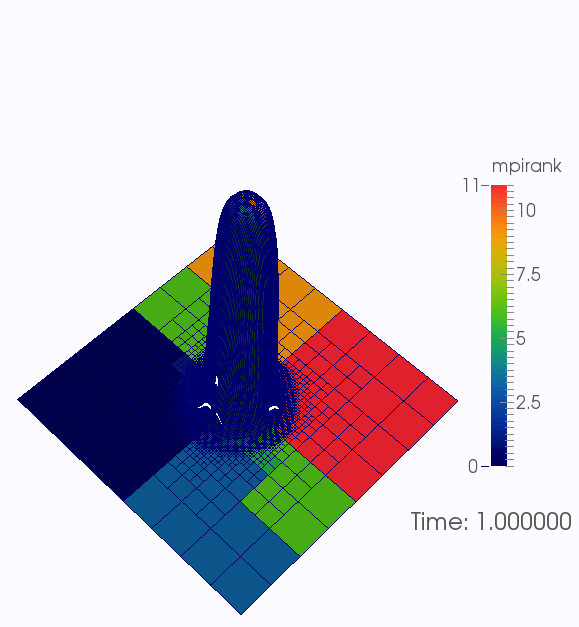}}
    \caption{View of the adaptive meshing and domain decomposition, load-balancing and 2:1 balance for the disk advection test case. 
             $\Delta x_{max} = 2^{-3}$, $\Delta x_{min} = 2^{-8}$. 2nd-order MUSCL-Hancock scheme.}
    \label{mesh-disk}
  \end{center}
\end{figure}

We consider again the transport problem at constant pressure and constant velocity of section~\ref{subsec-adv} with 
a sharp profile of volume fraction defined by 
$\alpha_0(\vect{x}) = 1 - \lambda$ if $|\vect{x} - \vect{x_0}| < 0.1 $, $\alpha_0(\vect{x})=\lambda$ otherwise.
The domain is periodic.  AMR is governed by the 
\textbf{$\rho$-gradient} criterion with the same refinement threshold as in the 
\ref{subsec-adv} case. 
Figure~\ref{mesh-disk} shows the resulting profiles obtained with the MUSCL-Hancock scheme
at several instants with a color representation of the 12 MPI processes domain decomposition.
The refinement criterion and the 2:1 balance property are well-managed by \pf{}. 

\subsubsection{Adapted versus uniform meshes}
In this section we compare results obtained with uniform grids and adapted meshes in order to assess the ability of the AMR procedure to act as a
compression technique, preserving accuracy while decreasing the computational needs. 
In the following, $\dxmin$ is the space step of the reference uniform mesh. 
It is equal to the size of the most refined cell in the AMR simulation and it is fixed for a series of simulation. 
On the other hand $\dxmax$ is the largest allowed space step for the AMR simulation. It varies from $\dxmin$ to $2^6 \dxmin$, so that 
the so called \textit{level of compression}, $\log_2\left(\dxmax\right)-\log_2\left(\dxmin\right)$ varies from 0 to 6.
The second order scheme and the \textbf{$\rho$-gradient} refinement criterion have been used to perform the different simulations.

\begin{figure}
  \begin{center}
        \subfloat[$\xi = 5 \times10^{-4}$]{\includegraphics[width=8cm]{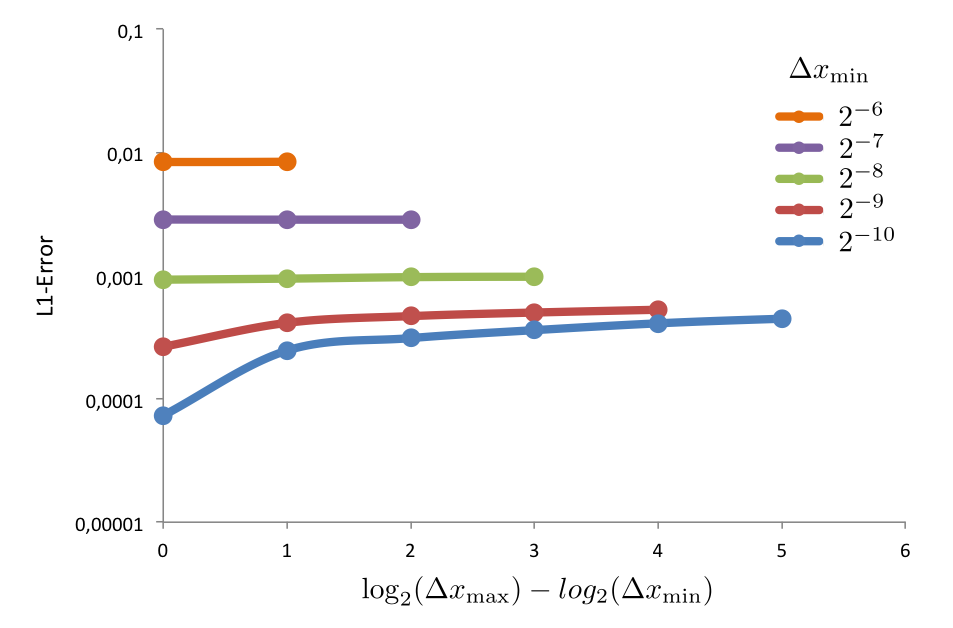}
    \label{compare-error-e4}}
    \hfill
    \subfloat[$\xi = 5 \times10^{-5}$]{\includegraphics[width=8cm]{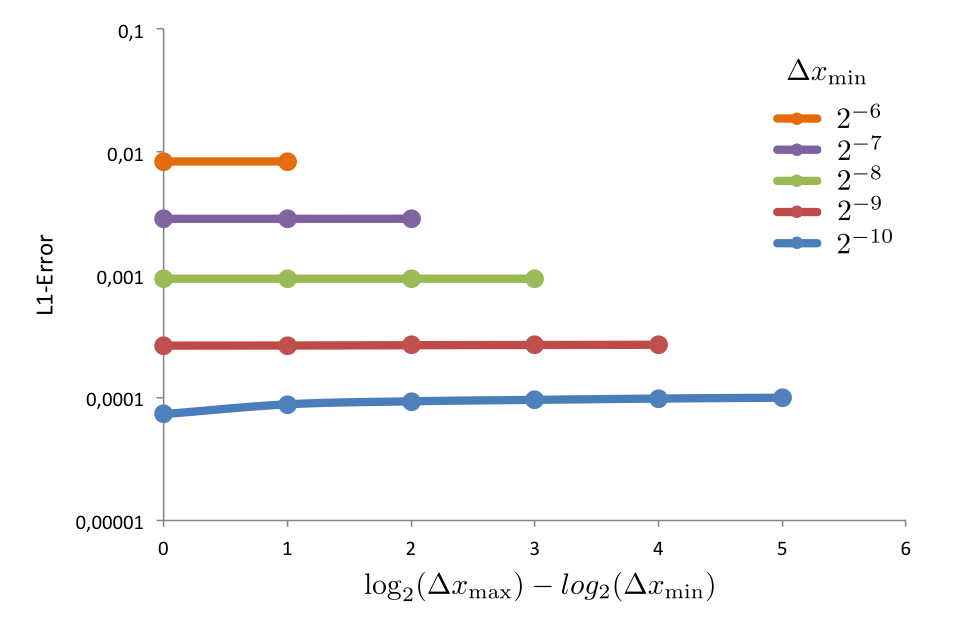}
    \label{compare-error-e5}}
    \caption{L$^1$-error versus level of compression of the mesh. Each compressed mesh is compared with its equivalent uniform mesh ($\log_2\left(\dxmax\right)-\log_2\left(\dxmin\right) = 0$) given in the same color.
             The study is lead for two different values of the threshold $\xi$ for \textbf{$\rho$-gradient} 
             refinement criterion.}
    \label{fig-criteria}
  \end{center}
\end{figure}

Figure~\ref{fig-criteria} shows the evolution of the L$^1$-error with the \textit{level of compression}
for two different values of the refinement threshold $\xi$. 
There, we see that for too large a refinement criterion, 
the compression error may prevail over the scheme consistency error when the space step $\dxmin$ goes to zero. 
This can be seen in Figure \ref{compare-error-e4} for $\dxmin = 2^{-9}$ and $\dxmin = 2^{-10}$ (red and blue lines),
where the L$^1$ errors of the compressed meshes are significantly higher than the L$^1$ error of the equivalent 
uniform mesh and do not seem to decrease with finer $\dxmin$.
However, decreasing the threshold $\xi$ for the refinement criterion enables to recover the 
expected accuracy of the \textit{compressed simulations}, as shown in Figure \ref{compare-error-e5}.
Then, there exists a subtle equilibrium for the refinement criterion: too small a value implies refinement everywhere and 
cancelation of the advantages of the AMR technique, whereas too large a value makes the compression error so large that 
mesh convergence is lost.

\begin{figure}
  \begin{center}
    \subfloat[At a fixed highest level of refinement $\dxmin$, we compare the total CPU time (full line and marks) 
      and compression rates (dashed lines) of the uniform mesh (in black) and the meshes with 1, 2 and 6 levels of compression
      ($\dxmax = 2^{1,2\text{ or } 6}\dxmin$). 
      Performed on 1 MPI process.]{\includegraphics[width=8cm]{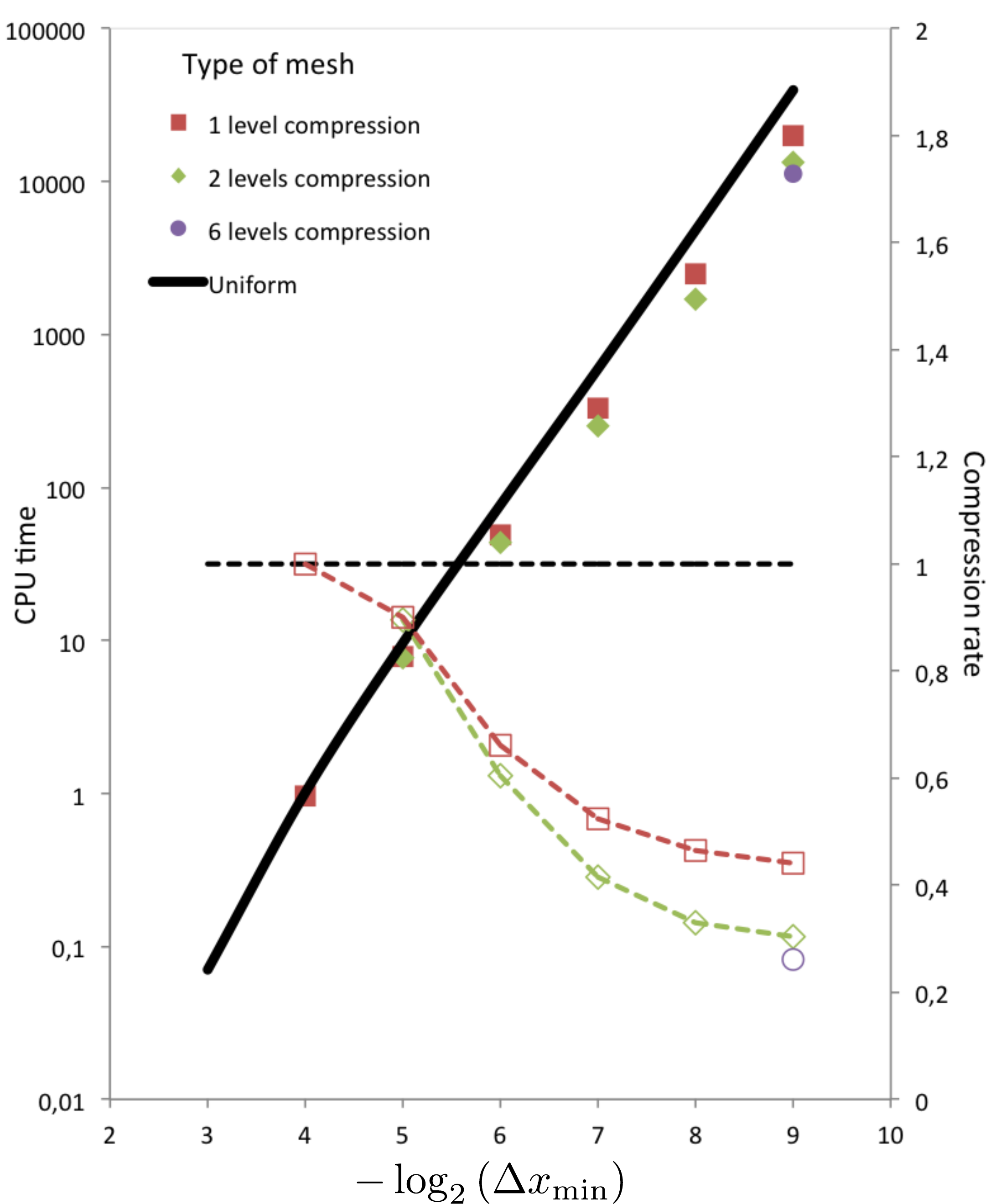}\label{fig-CPUtime}}
    \hfill
        \subfloat[Mean computation times (total time, time spent in finite volume solver and time for mesh management) per quadrant, for different problem sizes but with a same amount of work per process (weak scaling). ]{\includegraphics[width=8cm]{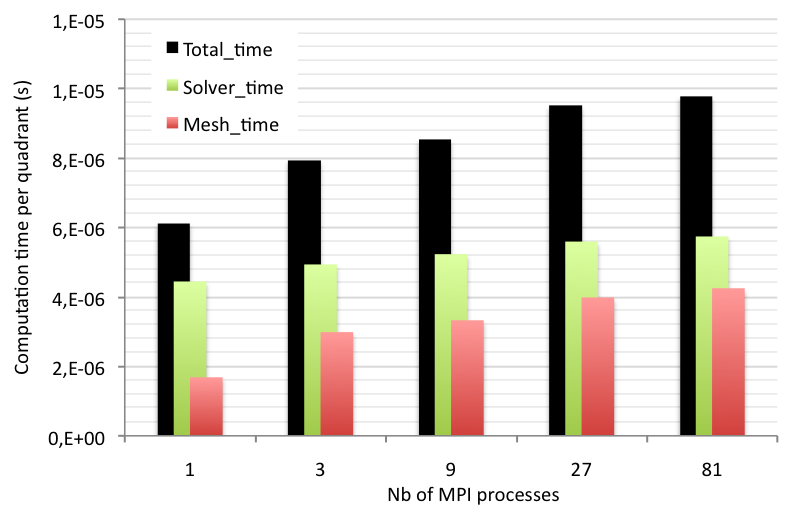}\label{weak-scaling}}
    \caption{
      Adaptive mesh resolution and computational costs}
  \end{center}
\end{figure}

When the refinement criterion is sufficiently small, the accuracies of the uniform and AMR solutions are comparable,
and the computational time on the compressed mesh is indeed better.
This is illustrated in Figure~\ref{fig-CPUtime}. This figure also displays the resulting compression rate, 
namely the ratio between the number of cells in the compressed mesh and in the equivalent uniform mesh. 
However, due to the rather large level of diffusion in the $\alpha$-advection test case, the compression rate is not as high as 
expected. We think that a less diffusive numerical scheme would enlarge uniform regions and therefore improve the AMR efficiency, 
needed for example in the case of the dynamics of a sharp interface between two phases.
Anyway, even though the AMR technique brings a certain overload for the management of 
the non-uniform mesh, the high compression rate allows an overall gain 
in terms of CPU time. 

\subsubsection{Parallel performance}
In this section, we present few cases testing some aspects of the parallel performance. 
Nonetheless, we need to emphasize that 
the code here is a raw first version that did not benefit from any optimization. 
It may be significantly improved in term of computational efficiency.

\espace

\textbf{Strong scaling}

\begin{figure}
  \begin{center}
    \subfloat[Strong scaling of total time of computation, of solver resolution and of mesh adaptation algorithms for about $4.6\times10^6$ cells.]{\includegraphics[width=7cm]{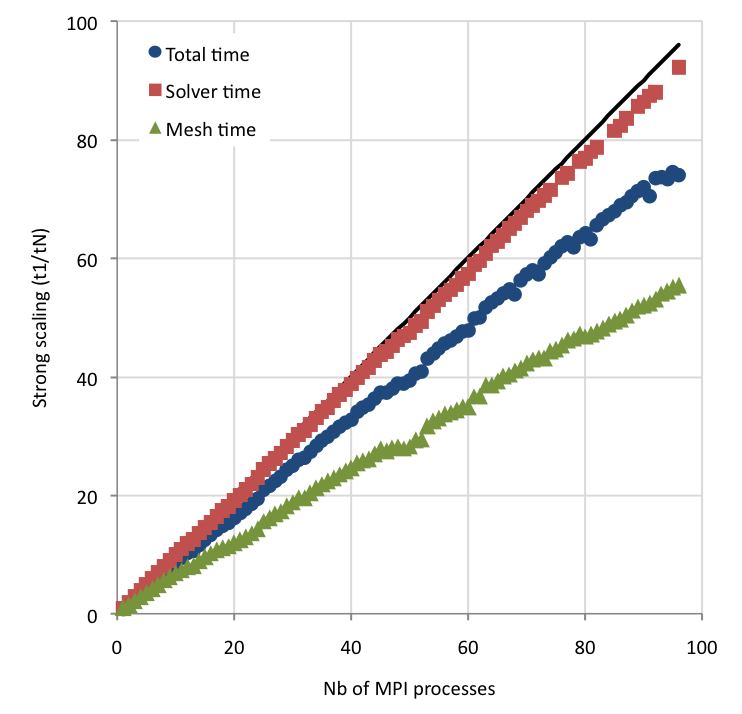} \label{strong-scaling}}
    \hfill
    \subfloat[Repartition of the computational time among the main tasks of the code according to meshes of different sizes and compression rates. 24 MPI-processes were used.]{\includegraphics[width=9cm]{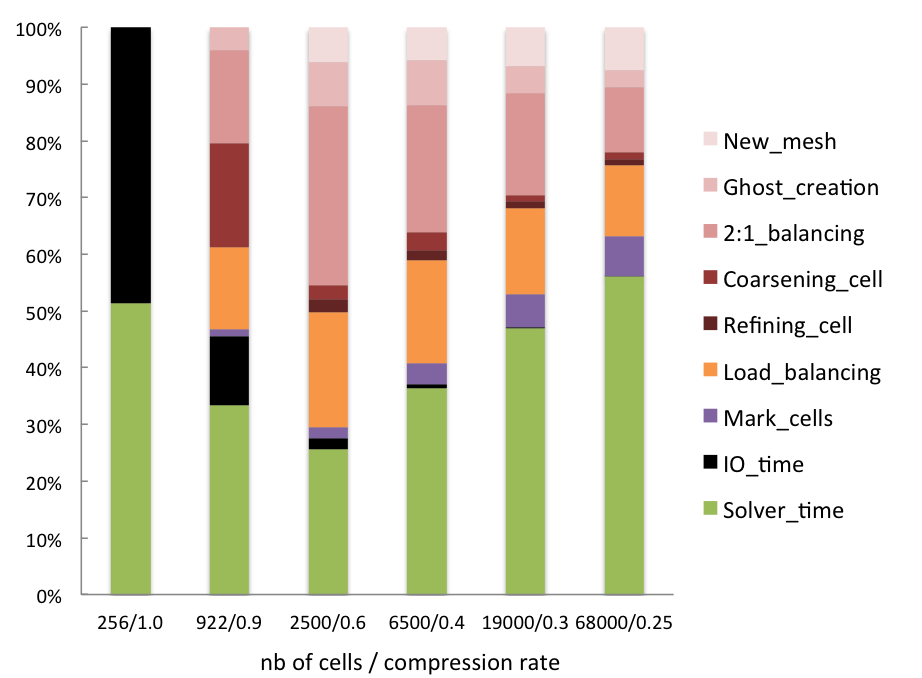}    \label{cost-of-work}}
    \caption{Strong scaling and time repartition for the $\alpha$-advection simulation.}
  \end{center}
\end{figure}

We use our $\alpha$-profile transport test
with a number of MPI processes ranging from 1 to 96. The runs are set in order to preserve the total number of cells approximately equal to $4.6\times10^6$ cells. Figure~\ref{strong-scaling} allows to evaluate the resulting speed-ups: for a low number of MPI processes the speed-up is very close to 1. However, in our case, for greater numbers of processes, the number of cells handled by the solver in each process is not sufficient to match the communication cost that becomes predominant. Indeed, as shown in Figure~\ref{mesh-disk}, 
the domain decomposition by even split of the z-curve does not always provide convex domains (some are not even connected). 
Therefore, the more MPI processes involved, the more the domain space is fragmented 
and the ratio of cells at the frontiers of each subdomain by its total number of cells increases, 
and so does the communications between the subdomains.

\espace

\textbf{Basic code profiling}

We perform an elementary profiling analysis in order to compare the CPU time allocated to the 
adaptation process versus the time spent in the Finite Volume solver 
given meshes that are successively refined, thus increasing their number of cells,
for a fixed number of 24 MPI processes. 
In Figure~\ref{cost-of-work} we display different tasks identified in the code.
We can see that the part dedicated to the mesh management given by the seven first colors (from light pink to dark red) 
decreases when the number of cells in the mesh increases. 
In particular, the part dedicated to the 2:1 balancing task, the longest one, is significantly reduced.

To sum up, Figures~\ref{strong-scaling} and \ref{cost-of-work} show that \pf{}~has good computation efficiency 
for important enough work loads, \textit{i.e.} for a high number of cells in the mesh. 
When the work load of each process is too low, 
an excessive time seems to be spent in communications between the processes compared to the time spent in the solver.

\espace

\textbf{Weak scaling}

We now evaluate the evolution of the computational time when increasing the number of working processes 
at a constant workload (see Figure~\ref{weak-scaling}). 
We maintain a number around $1.2\times10^4$ cells (the number of cells changes during the computation 
due to diffusion) managed by each process by increasing the global number of cells. The times of 
computation per quadrant are averaged over the 1000 first time steps.
While we did not succeed in preserving an exactly constant computational time per 
quadrant, 
the results are good and agree with similar results already obtained in \cite{burstedde11}.

%===========================================================
\subsection{2D and 3D gravity driven two-phase flows}
\label{par:3d}
In this section, we take the gravity source term $\S = (0,0,0,-\rho g, 0)^T$ into account by adding an additional operator 
$\Box_S^{\Delta t/2}$ in our splitting sequence \eqref{eq-update} or \eqref{eq-updateStrang}, following standard lines. 
For a given discrete state $\widetilde{\W}_i$ We set $\Box_S^{\Delta t/2} \widetilde{\W}_i 
= 
[
\widetilde{\rho}_i, 
\widetilde{(\rho Y)}_i, 
\widetilde{(\rho u_x)}_i, 
\widetilde{(\rho u_y)}_i - \widetilde{\rho}_i g\Delta t / 2,
\widetilde{(\rho u_z)}_i
]
$.
For three-dimensional problems, the overall splitting strategy becomes:
$$ 
\Box_X^{\Delta t/2}
\Box_Y^{\Delta t/2}
\Box_S^{\Delta t/2}
\Box_Z^{\Delta t/2}
\Box_Z^{\Delta t/2}
\Box_Y^{\Delta t/2}
\Box_S^{\Delta t/2}
\Box_X^{\Delta t/2}.
$$
We emphasize that the following tests aim at assessing the overall behavior of the code. 
While the physical behavior of the solutions is roughly correct, a more careful setup and systematic comparison with physical observations are still to be conducted in order to obtain a thorough validation.

\espace

\textbf{2D bubble drop test}

We consider the simulation of a falling drop of liquid (fluid~2) surrounded by a gas (fluid~1) 
toward a resting free surface separating a liquid bath from the gas. 
At $t=0$, we suppose that $\rho_1 = 1.0 (kg/m^3)$ and $\rho_2 = 1.0 \times10^3 (kg/m^3)$. 
We use solid wall boundary conditions.

We choose to use 
the \textbf{$\alpha$-gradient} refinement criterion, with $\xi = 5 \times10^{-4}$ 
in order to refine the mesh mainly in the vicinity of the gas/liquid interface.
\begin{figure}
\begin{center}
	\subfloat{\includegraphics[width=.32\textwidth]{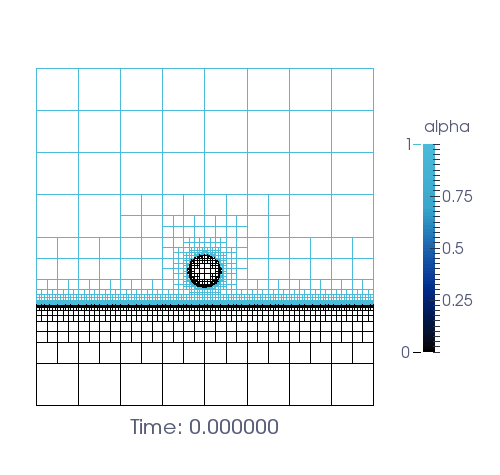}}
	\subfloat{\includegraphics[width=.32\textwidth]{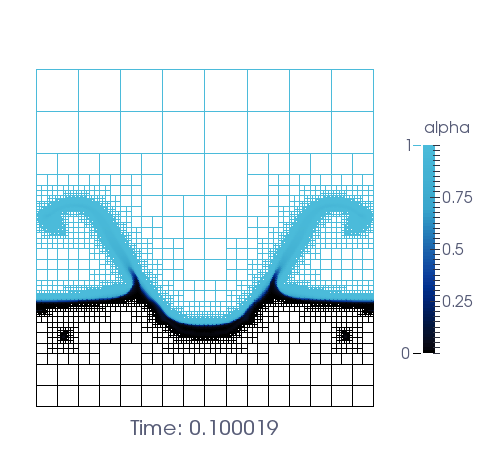}}
	\subfloat{\includegraphics[width=.32\textwidth]{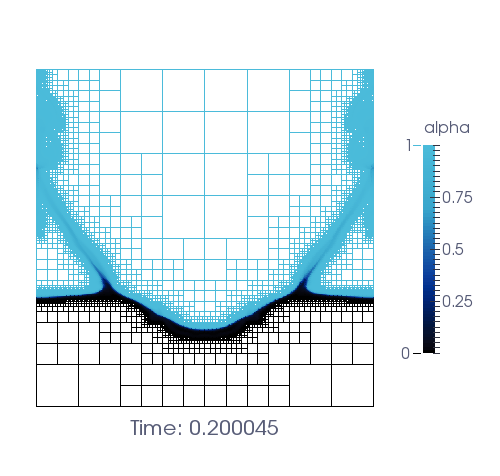}}\\
	\subfloat{\includegraphics[width=.32\textwidth]{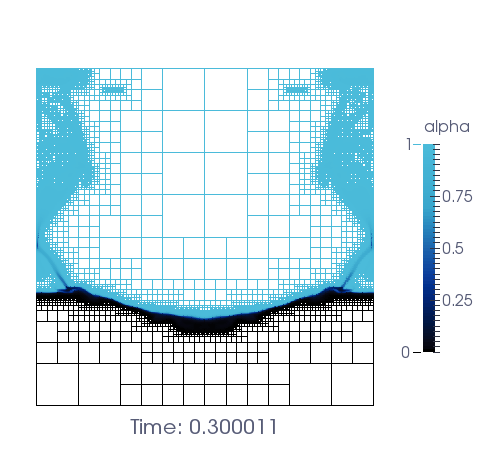}}
	\subfloat{\includegraphics[width=.32\textwidth]{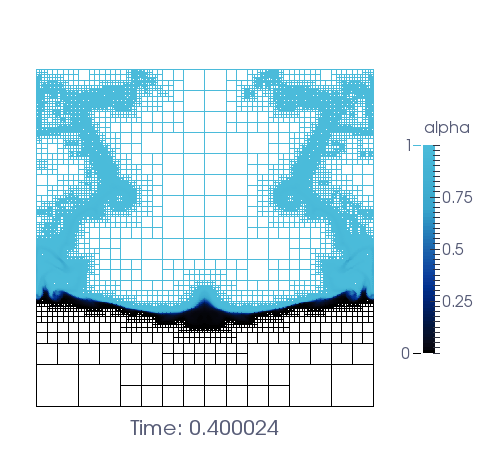}}
	\subfloat{\includegraphics[width=.32\textwidth]{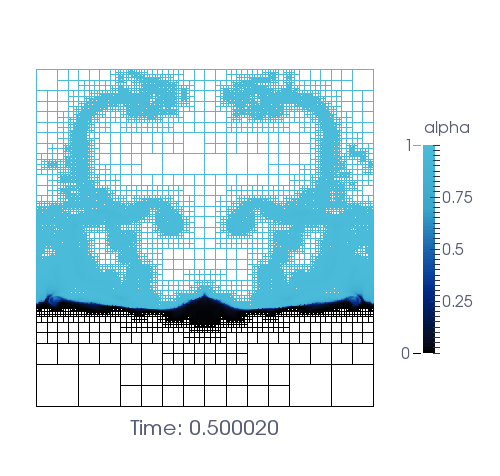}}\\
\caption{Simulation of a liquid drop falling onto a free surface with the refinement criterion \textbf{$\alpha$-gradient}. Mapping of the volume fraction $\alpha$. }
\label{fig-drop}
\end{center}
\end{figure}
The simulation is performed on a mesh with a minimum refining level of 3 and maximum of 9. 
Using 48 MPI processes, the time of computation is about 15 minutes on the computing 
Mesocenter of Ecole Centrale Paris, which is an Altix ICE 8400 LX. Each node is composed of two six-core Intel Xeon X5650 processors.
Figure~\ref{fig-drop} provides the mapping of the computed $\alpha$ at several instants along with the
adapted mesh. The mesh hits the finest refinement level in the vicinity of the interface where gradients of $\alpha$ are strong. 
Due to the numerical diffusion, we see that far from the interface gradients of $\alpha$ are detected and the mesh is also refined. 

\espace

\textbf{3D dam break test}

\begin{figure}
\begin{center}
	\subfloat{\includegraphics[width=8cm]{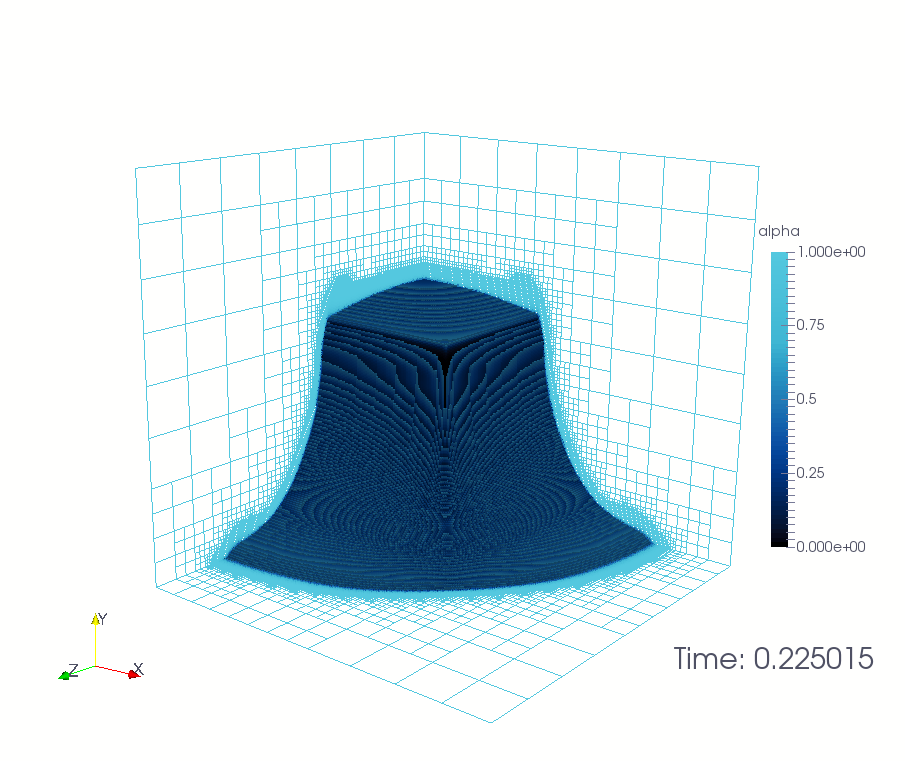}}
	\subfloat{\includegraphics[width=8cm]{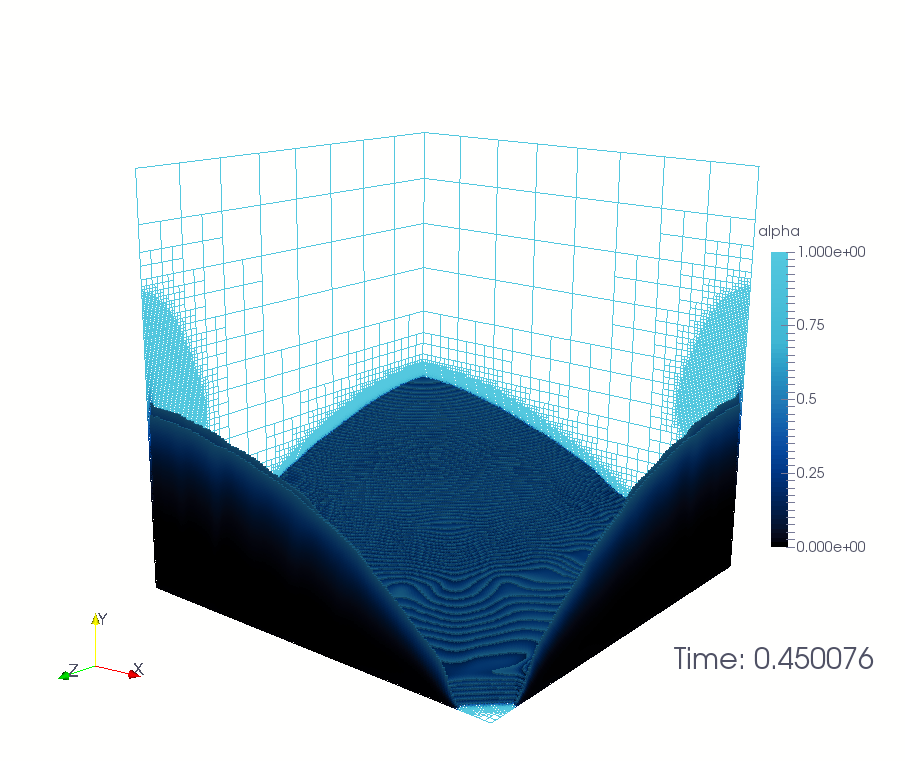}}\\
	\subfloat{\includegraphics[width=8cm]{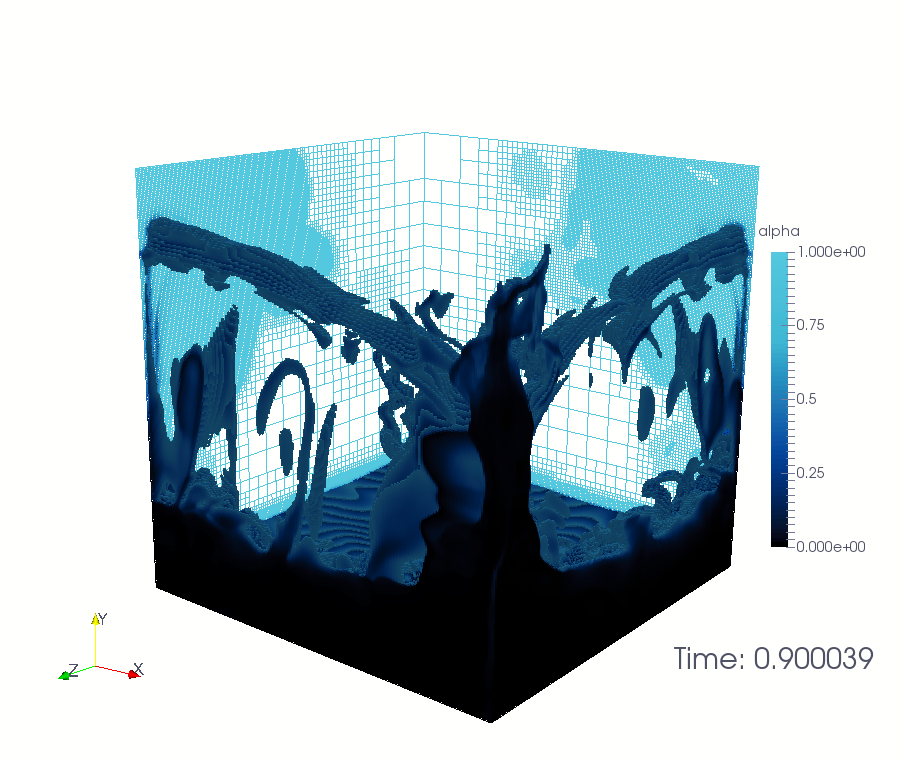}}
	\subfloat{\includegraphics[width=8cm]{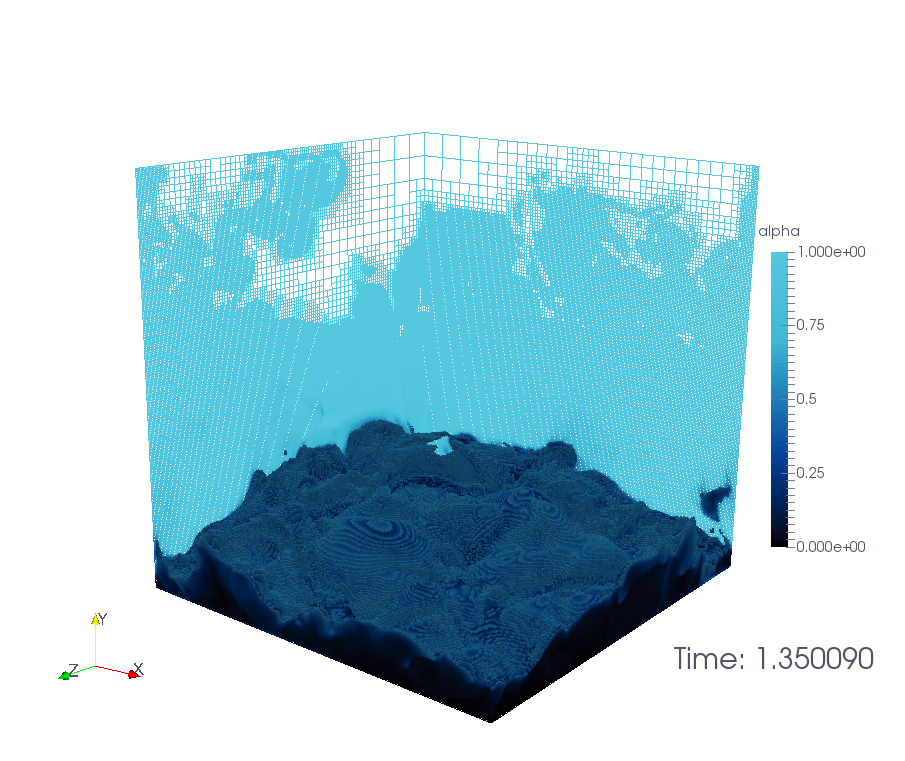}}
\caption{Simulation of a dam break. View of mesh and volume fraction $\alpha < 0.5$. Refinement criterion is \textbf{$\alpha$-gradient}.}
\label{fig-dam1}
\end{center}
\end{figure}

The second gravity-driven flow considered here deals with the evolution of a free surface in a dam break situation.
This problem has been studied in many works featuring simulations and experiments (see \textit{e.g.} \cite{Zhao11,Champmartin14}).
We show in Figure \ref{fig-dam1} the results of a 3D simulation. 
The computation was performed on 64 nodes of 8 CPU cores each of the computing Mesocenter of 
Centrale Paris.
A physical time of 1.5s has been reached: typically, given our initial conditions, the flow of liquid reaches
the opposite side of the domain within 0.3s and a second wave comes back within 1s.
The whole computation took about 4h for $5.59 \times10^5$ 
iterations. The highest level of refinement is 8 and the lowest 3. At the beginning of the 
simulation, the number of cells was $5.88 \times10^4$ ($3.5 \times 10^{-3}$ compression rate), 
at the end, due to diffusion and acoustic effects, it 
reached $1.25 \times10^6$ ($7.4 \times 10^{-2}$ compression rate), while the equivalent uniform 
mesh would have around $1.7\times10^7$ cells.

The refinement criterion is still $\alpha$-gradient and leads to an affordable computational cost, 
whereas the resolution of the problem on the finest grid would require a much longer time as well 
as a much larger memory.
The volume fraction iso-surface $\alpha < 0.5$, standing for the liquid phase, as well as the 
mesh at the domain boundaries are represented in the subfigures 
of Figure \ref{fig-dam1} at several time steps.  
At time $t=0$, both fluids are still. Due to gravity, the liquid flows into the chamber.
These results are very encouraging, even if they require further validation, as already stated, and if the influence of the refinement criterion on the dynamics of the solution has to be studied carefully.

%/\/\/\/\/\/\/\/\/\/\/\/\/\/\/\/\/\/\/\/\/\/\/\/\/\/\/\/\/\/\/\/\/\/\/\/\/\/\/\/\/\/\/\/\/\/\/\/\/\/\/\/\/\/\/\/\/\/\/\/\/\/\/\/\/\/\/\/\/\
%===========================================================
%/\/\/\/\/\/\/\/\/\/\/\/\/\/\/\/\/\/\/\/\/\/\/\/\/\/\/\/\/\/\/\/\/\/\/\/\/\/\/\/\/\/\/\/\/\/\/\/\/\/\/\/\/\/\/\/\/\/\/\/\/\/\/\/\/\/\/\/\/\

\section{Conclusion and perspectives}

The AMR library \pf~brings solutions to some issues for tree-based AMR, ranging from 
mesh and data structure using linear arrays, to cache locality thanks to the interesting properties of
the \textit{z-order curve}, and parallel efficiency through load balancing.
Understanding the main functionalities of \pf~and testing its ease of use and basic performance were the main objectives of  the six weeks of CEMRACS 2014.  
Within this framework, we have achieved a first version of a code, using a finite volume scheme of the relaxation type, 
at first and second order in space and time, applied to a simple but representative two-fluid two-phase flow model.
The scheme has been verified through classical test cases (advection, shock tube, double rarefaction) and a convergence analysis has been conducted. 

Some very promising simulations in 2D and 3D have been achieved and the code possesses all the good features in terms of parallel efficiency and accuracy, 
which allow both conducting reasonable size computations within a short amount of time 
(typically on a Mesocenter type of machine where the AMR strategy and its implementation lead to a solution at the same level of accuracy as uniform meshes but with significant savings in 
computational cost and memory requirement), and envisioning large scale and efficient simulations on larger massively parallel machines.
Such conclusions can be drawn, even if the tool requires both further optimization  and 
detailed and thorough study in terms of validation and accuracy of refinement criteria for the two-phase test-cases under study.

Let us also underline that there are some issues, which were not tackled in the present study. Among them, the problems we have studied do not involve a very large spectrum of time scales in terms of the dynamics of the problem \cite{duarte:hal-00727442} and the issue of local time stepping/multi-scale treatment will require some effort. Higher order numerical method will also require adapting the strategy proposed in the paper, as well as solving for elliptic equations such as in plasma physics and low speed flows (Low Mach approximation or incompressible flows). Such issues, even if interesting, were out of reach during the time of the project.

%%%%%%%%%%%%%%%%%%%%%%%%%%%%%%%%%%%%%%%%%%%
%%%%%%%%%%%%%%%%%%%%%%%%%%%%%%%%%%%%%%%%%%%

\subsection*{Acknowledgement}
The support of EM2C laboratory and of Maison de la Simulation for the CEMRACS project are gratefully acknowledged. The Ph.D. 
 of F. Drui is funded by a CEA/DGA (Direction Générale de l'Armement - French Department of Defense) grant.  
The use of the computational Mesocenter of Ecole Centrale Paris for some of the simulations is also gratefully acknowledged.

\bibliographystyle{plain}
\bibliography{biblio}

\begin{thebibliography}{10}

\bibitem{chombodesign}
M.~Adams, P.~Colella, D.~T. Graves, J.N. Johnson, N.D. Keen, T.~J. Ligocki. D.
  F. Martin.~P.W. McCorquodale, D.~Modiano.~P.O. Schwartz, T.D. Sternberg, and
  B.~Van Straalen.
\newblock Chombo software package for amr applications - design document.
\newblock {\em Lawrence Berkeley National Laboratory Technical Report}, 2013.

\bibitem{koumoutsakos2006}
M.~Bergdorf and P.~Koumoutsakos.
\newblock A {L}agrangian particle-wavelet method.
\newblock {\em Multiscale Model. Simul.}, 5(3):980--995, 2006.

\bibitem{berger89}
M.~J. Berger and P.~Collela.
\newblock Local adaptive mesh refinement for shock hydrodynamics.
\newblock {\em Journal of Computational Physics}, 82(1):64--84, 1989.

\bibitem{berger84}
M.~J. Berger and J.~Oliger.
\newblock Adaptive mesh refinement for hyperbolic partial differential
  equations.
\newblock {\em Journal of Computational Physics}, 53(3):484--512, 1984.

\bibitem{Champmartin14}
A.~Bernard-Champmartin and F.~De Vuyst.
\newblock A low diffusive lagrange-remap scheme for the simulation of violent
  air–water free-surface flows.
\newblock {\em Journal of Computational Physics}, 274(0):19 -- 49, 2014.

\bibitem{berthon06}
C.~Berthon.
\newblock Why the muscl-hancock scheme is {L}$^1$-stable.
\newblock {\em Numerische Mathematik}, 104:27--46, 2006.

\bibitem{bouchut}
F.~Bouchut.
\newblock {\em Nonlinear Stability of Finite Volume Methods for Hyperbolic
  Conservation Laws and Well-Balanced Schemes for Sources}.
\newblock Birkh{\"a}user, 2004.

\bibitem{brix2009parallelisation}
K.~Brix, S.~S. Melian, S.~M{\"u}ller, and G.~Schieffer.
\newblock Parallelisation of multiscale-based grid adaptation using
  space-filling curves.
\newblock In {\em ESAIM: Proceedings}, volume~29, pages 108--129. EDP Sciences,
  2009.

\bibitem{burstedde11}
C.~Burstedde, L.~C. Wilcox, and O.~Ghattas.
\newblock {\texttt{p4est}}: Scalable algorithms for parallel adaptive mesh
  refinement on forests of octrees.
\newblock {\em SIAM Journal on Scientific Computing}, 33(3):1103--1133, 2011.

\bibitem{caro}
F.~Caro, F.~Coquel, D.~Jamet, and S.~Kokh.
\newblock A simple finite-volume method for compressible isothermal two-phase
  flows simulation.
\newblock {\em International Journal of Finite Volume}, 3(1), 2006.

\bibitem{chalons2008}
C.~{Chalons} and J.-F. {Coulombel}.
\newblock {Relaxation approximation of the Euler equations.}
\newblock {\em {J. Math. Anal. Appl.}}, 348(2):872--893, 2008.

\bibitem{ONERA-chante}
G.~Chanteperdrix, P.~Villedieu, and J.P. Vila.
\newblock A compressible model for separated two-phase flows computations.
\newblock In {\em ASME Fluid Eng. Div. Summer Meeting 2002}, 2002.

\bibitem{koumoutsakos2007}
P.~Chatelain, G.-H. Cottet, and P.~Koumoutsakos.
\newblock Particle mesh hydrodynamics for astrophysics simulations.
\newblock {\em Internat. J. Modern Phys. C}, 18(4):610--618, 2007.

\bibitem{Cohen03}
A.~Cohen, S.~M. Kaber, S.~M\"{u}ller, and M.~Postel.
\newblock Fully adaptive multiresolution finite volume schemes for conservation
  laws.
\newblock {\em Mathematics of Computation}, 72:183--225, 2003.

\bibitem{postel09}
F.~Coquel, Q.L. Nguyen, M.~Postel, and Q.H. Tran.
\newblock Local time stepping applied to implicit-explicit methods for
  hyperbolic systems.
\newblock {\em Multiscale Model. Simul.}, 8(2):540--570, 2009/10.

\bibitem{koumoutsakoscottet2000}
G.-H. Cottet and P.~D. Koumoutsakos.
\newblock {\em Vortex methods}.
\newblock Cambridge University Press, Cambridge, 2000.
\newblock Theory and practice.

\bibitem{sechelles2015}
S.~Descombes, M.~Duarte, T.~Dumont, , T.~Guillet, V.~Louvet, and M.~Massot.
\newblock {Task-based adaptive resolution of time-space multi-scale
  reaction-diffusion systems on multi-core shared memory architectures}.
\newblock {\em SIAM Journal on Scientific Computing}, pages 1--24, 2015.
\newblock Submitted - Available on HAL
  https://hal.archives-ouvertes.fr/hal-01148617.

\bibitem{duarte:tel-00667857}
M.~Duarte.
\newblock {\em {Adaptive numerical methods in time and space for the simulation
  of multi-scale reaction fronts}}.
\newblock Th\`ese, {Ecole Centrale Paris}, December 2011.
\newblock \url{https://tel.archives-ouvertes.fr/tel-00667857}.

\bibitem{duarte:hal-00727442}
M.~Duarte, S.~Descombes, C.~Tenaud, S.~Candel, and M.~Massot.
\newblock {Time-space adaptive numerical methods for the simulation of
  combustion fronts}.
\newblock {\em {Combustion and Flame}}, 160(6):1083--1101, 2013.

\bibitem{duarte2012}
M.~Duarte, M.~Massot, S.~Descombes, C.~Tenaud, T.~Dumont, V.~Louvet, and
  F.~Laurent.
\newblock New resolution strategy for multiscale reaction waves using time
  operator splitting, space adaptive multiresolution, and dedicated high order
  implicit/explicit time integrators.
\newblock {\em SIAM J. Sci. Comput.}, 34(1):A76--A104, 2012.

\bibitem{Dubey20143217}
A.~Dubey, A.~Almgren, J.~Bell, M.~Berzins, S.~Brandt, G.~Bryan, P.~Colella,
  D.~Graves, M.~Lijewski, F.~Löffler, B.~O’Shea, E.~Schnetter, B.~Van
  Straalen, and K.~Weide.
\newblock A survey of high level frameworks in block-structured adaptive mesh
  refinement packages.
\newblock {\em Journal of Parallel and Distributed Computing}, 74(12):3217 --
  3227, 2014.
\newblock Domain-Specific Languages and High-Level Frameworks for
  High-Performance Computing.

\bibitem{harten95}
A.~Harten.
\newblock Multiresolution algorithms for the numerical solution of hyperbolic
  conservation laws.
\newblock {\em Comm. Pure and Applied Math.}, 48:1305--1342, 1995.

\bibitem{burstedde12}
T.~Isaac, C.~Burstedde, and O.~Ghattas.
\newblock Low-cost parallel algorithms for 2:1 octree balance.
\newblock In {\em Parallel Distributed Processing Symposium (IPDPS), 2012 IEEE
  26th International}, pages 426--437, May 2012.

\bibitem{ji}
H.~Ji, F.-S. Lien, and E.~Yee.
\newblock A new adaptive mesh refinement data structure with an application to
  detonation.
\newblock {\em Journal of Computational Physics}, 229(23):8981--8993, November
  2010.

\bibitem{khokhlov}
A.~M. Khokhlov.
\newblock Fully threaded tree for adaptive refinement fluid dynamics
  simulations.
\newblock {\em Journal of Computational Physics}, 143(2):519--543, July 1998.

\bibitem{koumoutsakos2005}
P.~Koumoutsakos.
\newblock Multiscale flow simulations using particles.
\newblock In {\em Annual review of fluid mechanics. {V}ol. 37}, volume~37 of
  {\em Annu. Rev. Fluid Mech.}, pages 457--487. Annual Reviews, Palo Alto, CA,
  2005.

\bibitem{leveque}
R.~J. Leveque.
\newblock {\em Finite-Volume Methods for Hyperbolic Problems}.
\newblock Cambridge University Press, 2004.

\bibitem{monaghan}
J.~J. Monaghan.
\newblock Smoothed particle hydrodynamics.
\newblock {\em Annual Review of Astronomy and Astrophysics}, 30:543--574, 1992.

\bibitem{muller}
S.~M{\"u}ller.
\newblock {\em Adaptive Multiscale Schemes for Conservation Laws}, volume~27.
\newblock Springer, 2003.
\newblock {Ed. T. J. Barth, M. Griebel, D. E. Keyes, R. M. Nieminen , D. Roose
  and T. Schlick}.

\bibitem{muller07}
S.~M{\"u}ller and Y.~Stiriba.
\newblock Fully adaptive multiscale schemes for conservation laws employing
  locally varying time stepping.
\newblock {\em J. Sci. Comput.}, 30(3):493--531, 2007.

\bibitem{popinet03}
S.~Popinet.
\newblock Gerris: a tree-based adaptive solver for the incompressible euler
  equations in complex geometries.
\newblock {\em Journal of Computational Physics}, 190(2):572--600, September
  2003.

\bibitem{koumoutsakos2011}
D.~Rossinelli, B.~Hejazialhosseini, D.~G. Spampinato, and P.~Koumoutsakos.
\newblock Multicore/multi-{GPU} accelerated simulations of multiphase
  compressible flows using wavelet adapted grids.
\newblock {\em SIAM J. Sci. Comput.}, 33(2):512--540, 2011.

\bibitem{suliciu}
I.~Suliciu.
\newblock On modelling phase transitions by means of rate-type constitutive
  equations, schock wave structure.
\newblock {\em International Journal of Engineering Science}, 1:829--841, 1990.

\bibitem{teyssier02}
R.~Teyssier.
\newblock Cosmological hydrodynamics with adaptive mesh refinement. a nex high
  resolution code called ramses.
\newblock {\em Astronomy and Astrophysics}, 385:337--364, 2002.

\bibitem{toro}
E.~F. Toro.
\newblock {\em Riemann Solvers and Numerical Methods for Fluid Dynamics - A
  Practical Introduction}.
\newblock Springer, 3rd edition, 2009.

\bibitem{toro1994}
E.F. {Toro}, M.~{Spruce}, and W.~{Speares}.
\newblock {Restoration of the contact surface in the HLL-Riemann solver.}
\newblock {\em {Shock Waves}}, 4(1):25--34, 1994.

\bibitem{leer84}
B.~van Leer.
\newblock On the relation between the upwind-differencing schemes of godunov,
  engquist-osher and roe.
\newblock {\em SIAM Journal on Scientific and Statistical Computing},
  5(1):1--20, 1984.

\bibitem{wood}
A.~B. Wood.
\newblock {\em A textbook of Sound}.
\newblock The Macmillan Company, 1930.

\bibitem{Zhao11}
X.-Z. Zhao.
\newblock Validation of a cip-based tank for numerical simulation of free
  surface flows.
\newblock {\em Acta Mechanica Sinica}, 27(6):877--890, 2011.

\end{thebibliography}

\end{document}